\documentclass[10pt]{article}

\usepackage{amssymb}
\usepackage{amsmath}
\usepackage{theorem}
\usepackage[dvips]{graphicx}
\usepackage{epstopdf}
\usepackage{eepic}
\usepackage{enumitem}
\usepackage{diagbox}

\newtheorem{theorem}{Theorem}[section]
\newtheorem{proposition}[theorem]{Proposition}
\newtheorem{lemma}[theorem]{Lemma}

\newtheorem{definition}[theorem]{Definition}

{\theorembodyfont{\rmfamily}\newtheorem{example}[theorem]{Example}}
{\theorembodyfont{\rmfamily}\newtheorem{remark}[theorem]{Remark}}
\newenvironment{keywords}{ {\bf Key words.} }{}

\newcommand{\Int}{\int\limits}

\title{Robin-Dirichlet alternating iterative procedure for solving  the Cauchy problem for
 Helmholtz equation in an unbounded domain}
 
\author{Pauline Achieng$^{\rm a,b,}$\thanks{$^\ast$ Corresponding author. Email: achiengpauline2@gmail.com
$^{\rm b}${\em{Department of applied Mathematics,University of Nairobi, P.O. Box 30197,
Nairobi}} $^{\rm a}${\em{Department of Mathematics,  Link\"{o}ping University, SE-581 83 Link\"{o}ping, Sweden}}},
Fredrik Berntsson$^{\rm a}$ and Vladimir Kozlov$^{\rm a}$}
 
\date{\today}

\begin{document}

\thispagestyle{empty}
\maketitle

%
\begin{abstract}                                                     
We consider the Cauchy problem for the Helmholtz equation with a domain in $\mathbb R^d$, $d\geq 2$  with $N$ cylindrical outlets to infinity with  bounded inclusions in $\mathbb R^{d-1}.$ Cauchy data are prescribed on the boundary of the bounded domains and the aim is to find solution on the unbounded part of the boundary.
In 1989, Kozlov and Maz${'}$ya \cite{Kozlov:1989} proposed an alternating iterative method for solving Cauchy problems associated with elliptic, self-adjoint and positive-definite operators in bounded domains. Different variants of this method for solving Cauchy problems associated with Helmholtz-type operators exists. We consider the variant proposed by Mpinganzima et.al \cite{Lydie:2018} for bounded domains and  derive the necessary conditions  for the convergence of the  procedure  in  unbounded domains. For the numerical implementation, a finite difference method is used to solve the problem in a simple rectangular domain in $\mathbb R^2$ that represent  a truncated infinite strip. The numerical results shows that by appropriate truncation of the domain and with appropriate choice of the Robin parameters $\mu_0$ and $\mu_1$, the Robin-Dirichlet alternating iterative procedure is convergent.
\end{abstract}


\section{Introduction}\label{sec:introduction}
Let $\Omega$ be a domain in $\mathbb R^d$, $d\geq 2$, with $C^2$ boundary and with $N$ cylindrical outlets to infinity, i.e. for sufficiently large $|x|$ the domain $\Omega$ coincides with the union of $N$ disjoint cylinders ${\mathcal C}^{(j)}$, $j=1,\ldots,N$, which can be described in a certain cartesian coordinates $x^{(j)}=(y^{(j)},z^{(j)})$,  as
$$
{\mathcal C}^{(j)}=\{x^{(j)}\,:\,y^{(j)}\in\omega^{(j)},\,z^{(j)}\in\mathbb R\},
$$
where the cross-sections $\omega^{(j)}$ are bounded domains in $\mathbb R^{d-1}$ with $C^2$ boundaries. We denote the boundary of $\Omega$ by $\Gamma$. We assume that a certain bounded\footnote{This is a set where measurements are taken and it is reasonable to assume it bounded} open set $\Gamma_0$ is chosen on the boundary $\Gamma$ and the boundary of this set is of class $C^2$ also. Let also $\Gamma_1$ is the interior of $\Gamma\setminus\Gamma_0$.

We consider the following Cauchy problem for the Helmholtz equation
\begin{equation}\label{I1}
(\Delta+k^2)u=0\;\;\mbox{in $\Omega$}
\end{equation}
and
\begin{eqnarray}\label{I2}
&&u=f_0\;\;\mbox{on $\Gamma_0$},\nonumber\\
&&\partial_\nu u=g_0\;\;\mbox{on $\Gamma_0$},
\end{eqnarray}
where $k$ is a real scalar, $\nu$ is the outward unit normal to $\Gamma$, $\partial_\nu$ is the normal derivative and $(f_0,g_0)$ is a prescribed Cauchy data.

The Cauchy problem for the Helmholtz equation in bounded and unbounded domains  arises in many important physical applications, for instance in capacity problems or scattering of  acoustics or electromagnetic waves, see \cite{chen:1998, delillo:2001, delillo:2003, kress:1998, jones:1986}.

The Cauchy problem for the Helmholtz equation is an inverse problem and it is ill-posed. Small perturbation in the Cauchy data  $f_0$ and $g_0$ results into a big error in the solution and as a result classical numerical methods cannot be used to solve this problem. Regularization methods are instead used to solve inverse problems.

Over the years, much theoretical and numerical studies have been done on the Cauchy problem associated with the Helmholtz equation on bounded domains. These include both Tikhonov type regularization methods and iterative regularization methods. Lesnic et al. \cite{lesnic:2003b} and  Marin \cite{marin:2009a},    have solved the Cauchy problem associated with the Helmholtz equation using the conjugate gradient method (CGM) and the boundary element-minimal error method, respectively. Wei et al. \cite{wei:2009, qin2008tikhonov} solved the Cauchy problem associated with the Helmholtz-type equations by transforming the Cauchy problem into a moment problem and then applied  a Tikhonov type regularization method. Zhang et al. \cite{liu2017fourier} have solved the Cauchy problem for the Helmholtz equation using a Fourier-Bessel method.   Numerical methods for solving the Cauchy problem for two and three dimensional Helmholtz-type equations have also been studied by 
Lesnic et at.  \cite{marin2005method} and Marin  \cite{marin2005meshless} respectively. They used  the method of fundamental solutions (MFS) in conjunction with Tikhonov regularization method.

 Kozlov and Maz${'}$ya \cite{Kozlov:1989} developed the alternating iterative procedure for solving  linear elliptic partial differential equations. The alternating iterative procedure is applicable to equations where the operator is symmetric and positive in a certain sense.  The regularizing character is achieved by appropriately changing the boundary conditions and one advantage of this procedure is that it preserves the original operator.
 Kozlov et al. \cite{Kozlov:1991} used this procedure to solve the Cauchy problem for the Laplace equation and  the Lame${'}$  system.  Chapko et al. \cite{chapko2008alternating}  further  applied  the alternating iterative procedure to solve the Cauchy problem for the Laplace equation in a bounded domain with a cut. It has also been demonstrated that the alternating iterative procedure does not only work for linear elliptic partial differential equations but also for nonlinear elliptic partial differential equations, see  \cite{maxwell2008iterative, avdonin:2009iterative}. However, as mentioned above, the alternating iterative procedure converge if the operator is self-adjoint and  positive-definite. Example of  operators which do not fulfil this requirements are Helmholtz-type operators.  Marin et.al  \cite{marin2003alternating} used the alternating iterative procedure, numerically implemented using the  boundary element method (BEM) to solve the Cauchy problem for the Helmholtz equation with purely imaginary wavenumber, $k$ i.e they considered  the equation $( \Delta-k^2)u=0$ which is in fact a modified version of the Helmholtz equation. They noticed that the alternating iterative procedure applied to the Helmholtz equation does not always converge.
 Kozlov et al. \cite{johansson2009alternating}  modified the alternating iterative procedure to accommodate second order  elliptic operators which are self-adjoint but does not fulfil the condition of positivity.
 Mpinganzima  et al.  \cite{berntsson2014alternating, Lydie:2014b}  also presented other modifications of the alternating iterative procedure for Cauchy problem associated to the Helmholtz equation by an introduction of artificial boundary and boundary conditions which allow to treat all values of $k$ i.e the positivity condition introduced in   \cite{Kozlov:1989, Kozlov:1991}  is essentially selected. The latter,  \cite{Lydie:2014b},  involves employing  an operator equation formulation of the Robin-Dirichlet algorithm and using the conjugate gradient method in order to accelerate the slow convergence achieved in \cite{berntsson2014alternating}.
 In \cite{Lydie:2018},  Mpinganzima  et al. further presented a simpler modification of the alternating iterative procedure for Cauchy problem associated to the Helmholtz equation by replacing the Neumann-Dirichlet iterations by the Robin-Dirichlet iterations.  In \cite{achieng2020}, we presented an analysis of  Robin-Dirichlet alternating iterative procedure. We prove that the  Robin-Dirichlet alternating iterative procedure is in fact convergent for general elliptic operators provided that the parameters in the Robin conditions are appropriately chosen.  The precise behaviour of $k$  in the Helmholtz equation is also numerically investigated.

The aim of this paper is to derive the necessary conditions  for the convergence of the Robin-Dirichlet alternating iterative procedure for solving  the Cauchy problem for the Helmholtz equation in  unbounded domains.
 In unbounded domains, for example in the cylinder ${\mathcal C}^{(j)}$ considered in our problem, the continuous 
 spectrum of the Dirichlet-Laplacian in  $\Omega$  coincides with  $[\min_j\lambda^{(j)}_0,\infty)$ where $\lambda_0^{(j)}$ is the first eigenvalue of the Dirichlet-Laplacian in the  cross-section $\omega^{(j)}$,  see \cite{jones1953eigenvalues, nazarov2010variational}. If there is no discrete spectrum below   $\min_j\lambda^{(j)}_0$, then we prove the convergence of the Robin-Dirichlet alternating iterative procedure if $k^2<\lambda^{(j)}_0$  for all j. However, if there are eigenvalues below $\min_j\lambda^{(j)}_0$ then we prove
 the convergence of the Robin-Dirichlet alternating iterative procedure if $k^2< \Lambda_0$ where $ \Lambda_0$ is the smallest eigenvalue in the discrete spectrum. The convergence analysis of the 
Robin-Dirichlet alternating iterative procedure is based on an analysis of the spectrum of the Laplacian operator in $\Omega$ with Dirichlet and Robin boundary conditions.

\subsection{Alternating iterative procedure}\label{sec:Two auxiliary problems and alternating iterative procedure}

As usual the notation $H^1(\Omega)$ corresponds to the Sobolev space of functions in $\Omega$ with finite norm
$$
||u||_{H^1(\Omega)}=\Big(\int_\Omega(|\nabla u|^2+|u|^2)dx\Big)^{1/2}.
$$
 Our main assumption concerning the parameter $k$ is the following: there exist a positive constant $\epsilon$ such that
\begin{equation}\label{I3}
\int_\Omega (|\nabla u|^2-k^2|u|^2)dx\geq \epsilon ||u||^2_{H^1(\Omega)}\;\;\mbox{for all $u\in H^1(\Omega,\Gamma)$},
\end{equation}
where $H^1(\Omega,\Gamma)$ is the subspace of function in $H^1(\Omega)$ vanishing on $\Gamma$.

In Lemma \ref{L2},  we give an equivalent version of the assumption (\ref{I3}). Namely, it is equivalent to the following: there exist positive constants $\mu_0$, $\mu_1$  and $\delta$ such that
\begin{eqnarray}\label{I4}
\int_\Omega (|\nabla u|^2-k^2|u|^2)dx+\mu_0\int_{\Gamma_0}|u|^2dS+\mu_1\int_{\Gamma_1}|u|^2dS\geq \delta ||u||^2_{H^1(\Omega)}
\end{eqnarray}
for all $u\in H^1(\Omega)$.

In order to describe the Robin-Dirichlet alternating iterative procedure, let us introduce two mixed  boundary value problems:
\begin{equation}\label{eq:auxiliary:2}
\begin{cases}
\Delta u + k^2u = 0 & \quad \mbox{in} \quad \Omega,\\
u = f   & \quad \mbox{on} \quad \Gamma_0,\\      
\partial_\nu u+ \mu_1 u =\eta  & \quad \mbox{on} \quad \Gamma_1,
\end{cases}
\end {equation}
and
\begin{equation}\label{eq:auxiliary:3}
\begin{cases}
\Delta u + k^2u = 0 & \quad \mbox{in} \quad \Omega,\\
  \partial_\nu u+  \mu_0 u= g      & \quad \mbox{on} \quad \Gamma_0,\\    
u = \phi    & \quad \mbox{on} \quad \Gamma_1,

\end{cases}
\end{equation}
Here, $f \in H^{\frac{1}{2}}(\Gamma_0) $, $\eta\in H^{-1/2 }(\Gamma_1) $, $g \in H^{-1/2 }(\Gamma_0) $ and  $ \phi \in H^{\frac{1}{2}}(\Gamma_1)$. From the assumption (\ref{I3}), or equivalently from (\ref{I4}), it follows well-posedness of these problems, see section 3.2.

The algorithm for solving (\ref{I1}), (\ref{I2}) is described as follows. We take $f = f_0$ and $g = g_0+ \mu_0 f_0$  where  $ f_0$ and $g_0$  are the Cauchy data given in (\ref{I2}) then;
\begin{itemize}
\item[(1)] The first approximation $u_{0}$ is obtained by solving (\ref{eq:auxiliary:2}) where $\eta\in H^{-1/2 }(\Gamma_1) $ is an arbitrary initial approximation of the Robin condition on $\Gamma_1$.
\item[(2)] Having constructed  $u_{2n}$, we obtain  $u_{2n+1}$ by solving (\ref{eq:auxiliary:3})  with $\phi=u_{2n}$ on $\Gamma_1$.
\item[(3)] We then obtain $u_{2n+2}$ by solving  (\ref{eq:auxiliary:2}) with $\eta =  \partial_\nu u_{2n+1}+  \mu_1 u_{2n+1}$
\end{itemize}
In section 3.3, we present a theorem on convergence of this algorithm. The proof basically follows the same lines as in the case of bounded domains, see \cite{achieng2020}.

 \section{About condition (\ref{I3})}
We denote by $\lambda_0^{(j)}$ the first eigenvalue of the operator  $-\Delta_{y^{(j)}}$ in the cross-section $\omega^{(j)}$  and $\lambda^{(j)}(\mu)$ the first eigenvalue of the operator $-\Delta_{y^{(j)}}$ in $\omega^{(j)}$ with the Robin boundary condition $(\partial_\nu+\mu)u=0$ on $\partial\omega^{(j)}$, $j=1,\ldots,N$.

As is known the continuous spectrum of $-\Delta$ in $\Omega$ lies in $[\min_j\lambda^{(j)}_0,\infty)$ and of the operator $-\Delta$ in $\Omega$ with the Robin boundary condition $(\partial_\nu+\mu)u=0$ on $\Gamma$ is located in $[\min_j\lambda^{(j)}(\mu),\infty)$. It can happen that there is also a discrete spectrum for both operators lying in $(-\infty,\min_j\lambda^{(j)}_0)$ and $(-\infty,\min_j\lambda^{(j)}(\mu))$ respectively.

\begin{lemma}\label{L1} The function $\lambda^{(j)}(\mu)$ is monotonically increasing with respect to $\mu$ and
\begin{equation}\label{M31a}
\lambda^{(j)}(\mu)\to \lambda^{(j)}_0\;\;\mbox{as $\mu\to\infty$},   \quad \mbox{for}  \;\; j=1,\ldots,N.
\end{equation}
\end{lemma}
{\bf Proof.} The eigenvalue $\lambda^{(j)}(\mu)$ of the operator $-\Delta_{y^{(j)}}$ in $\omega^{(j)}$ with the Robin boundary condition $(\partial_\nu+\mu)u=0$ on $\omega^{(j)}$  is given by
  \begin{equation}\label{eq:eigenpost:1}
\lambda^{(j)}(\mu)=\min_{\substack{u\in H^1(\Omega)\\   u\neq 0}} \frac{ \mu\int_{\partial \omega^{(j)}} u^2 dS +  \int_{ \omega^{(j)}}  \mid \nabla u \mid^2   \,dx}{ \int_{\omega^{(j)}} u^2 \,dx}.
\end{equation}
From (\ref{eq:eigenpost:1}), we see that the function $ \lambda^{(j)}(\mu)$ is non-negative and monotonically increasing with respect to $ \mu$. Then it remains to prove that $\lim_{ \mu \to\infty} \lambda^{(j)}(\mu)= \lambda^{(j)}_0$ which follows from Lemma 2.1 in \cite{berntsson2014alternating}.

\begin{lemma}\label{L0}
If condition { \rm(\ref{I3})} holds then 
\begin{equation}\label{I3a}
k^2<\lambda_0^{(j)},\;\;j=1,\ldots,N.
\end{equation}
\end{lemma}
{\bf Proof.}
Let $k^2\geq \lambda_0^{(j)}$ for a certain j.
We choose a test functions  in (\ref{I3}) in the form $u(x^{(j)})=\phi^{(j)}(y^{(j)})\eta(z^{(j)})$, where $\phi^{(j)}$ is an  eigenfunction of the Dirichlet-Laplacian in $ \omega^{(j)}$ corresponding  to $\lambda_0^{(j)}$ with the norm $||\phi^{(j)}||_{L_2(\omega^{(j)})}=1$. Substituting $u(x^{(j)})$ to the left-hand side of (\ref{I3}) we  obtain
\begin{eqnarray}\label{lem1a}
 \Int_{-\infty}^{\infty}  \Int_{\omega^{(j)}}\left( |\nabla (\phi^{(j)}\eta)|^2-k^2 (\phi^{(j)})^2(\eta)^2 \right) \,dy^{(j)}\,dz^{(j)}\nonumber\\
=  \Int_{-\infty}^{\infty} (\eta^\prime)^2+\left( \lambda_0^{(j)}-k^2\right)\eta^2 \,dz^{(j)}
\end{eqnarray}
Let the test function $\eta$ be equal to $1$ for $\tau+1<z^{(j)}<\tau+T-1$, $\eta(z^{(j)})=z^{(j)}-\tau$ for $\tau<z^{(j)}<\tau+1$, $\eta(z^{(j)})=-z^{(j)}+\tau+T$ for $\tau+T-1<-z^{(j)}<\tau+T$ and $0$ otherwise. Then
\begin{equation*}
 \Int_{-\infty}^{\infty} (\eta^\prime)^2 \,dz^{(j)}=2 < \frac{2}{T-2}\Int_{-\infty}^{\infty} \eta^2\,dz^{(j)}
\end{equation*}
and therefore the inequality (\ref{I3}) can not be true for all such test functions in the case $k^2\geq \lambda_0^{(j)}$. This prove this lemma.

%
  From Lemma 2.1 and inequality (\ref{I3a}), it follows that
\begin{equation}\label{I3aaa}
k^2<\lambda^{(j)}(\mu),\;\;j=1,\ldots,N,
\end{equation}
for large $\mu$.

Now we can prove the equivalence of the conditions (\ref{I3}) and (\ref{I4}).
\begin{lemma}\label{L2} Condition { \rm(\ref{I3})}  is equivalent to existence of positive constants $\mu$ and $\delta$ such that
\begin{equation}\label{I4a}
\int_\Omega (|\nabla u|^2-k^2|u|^2)dx+\mu\int_{\Gamma}|u|^2dS\geq \delta ||u||^2_{H^1(\Omega)}\;\;\mbox{for all $u\in H^1(\Omega)$}.
\end{equation}
\end{lemma}
{\bf Proof.} Clearly (\ref{I3}) follows from (\ref{I4a}). Let us prove the opposite implication.

We introduce
$$
\Lambda_0=\inf_{\substack{u\in H^1(\Omega,\Gamma)\\  {\Vert u \Vert_{L_2(\Omega)}=1}}}\int_\Omega |\nabla u|^2dx.
$$
Then
\begin{equation}\label{M31b}
\epsilon+k^2\leq \Lambda_0\leq \min_j\lambda^{(j)}_0.
\end{equation}
Furthermore, if we put
\begin{equation}\label{M31bb}
\Lambda(\mu)=\inf_{\substack{u\in H^1(\Omega)\\  {\Vert u \Vert_{L_2(\Omega)}=1}}}\Big(\int_\Omega |\nabla u|^2dx+\mu\int_\Gamma |u|^2dS\Big),
\end{equation}
where $\mu$ is non-negative, then on can see that
\begin{equation}\label{M31ba}
0\leq\Lambda(\mu)\leq \min_j\lambda^{(j)}(\mu)\leq\Lambda_0.
\end{equation}
Moreover $\Lambda(\mu)$ is a monotonically increasing function with respect to $\mu$. Then the required assertion will follow from
\begin{equation}\label{Jul11a}
\Lambda(\mu)\rightarrow \Lambda_0\;\;\mbox{as $\mu\to\infty$.}
\end{equation}
Let us prove (\ref{Jul11a}).

First, consider the case when there is a sequence $\{\mu_k\}_{k=1}^\infty$ such that $\mu_k\to\infty$ as $k\to\infty$ and $\Lambda(\mu_k)=\min_j\lambda^{(j)}(\mu_k)$. Since $\lambda^{(j)}(\mu)\to \lambda^{(j)}_0$ as $\mu\to\infty$, we get
$$
\min_j\lambda^{(j)}(\mu)\rightarrow\min_j\lambda^{(j)}_0\;\;\mbox{as $\mu\to\infty$.}
$$
Therefore $\Lambda(\mu_k)\rightarrow\min_j\lambda^{(j)}_0$ and we get
(\ref{Jul11a}) due to monotonicity of $\Lambda(\mu)$ and because of the right inequality in (\ref{M31b}) (in this case $\Lambda_0=\min_j\lambda^{(j)}_0$).

Second, suppose that
$\Lambda(\mu)<\min_j\lambda^{(j)}(\mu)$ for large $\mu$. Then $\Lambda(\mu)$ is an eigenvalue of the Laplacian in $\Omega$ with the Robin boundary condition $\partial_\nu u+\mu u=0$ on $\Gamma$. We denote by $u_\mu$ a corresponding eigenfunction normalized by $||u_\mu||_{L_2(\Omega)}=1$.
The function $u_\mu$ satisfies
$$
\Delta u_\mu=-\Lambda(\mu)u_\mu\;\;\mbox{in $\Omega$ and}\;\;\partial_\nu u_\mu+u_\mu=0\;\;\mbox{on $\Gamma$}
$$
and hence
\begin{equation}\label{Jul11b}
\int_{\Omega}|\nabla u_\mu|^2dx+\mu\int_\Gamma |u_\mu|^2dS=\Lambda_\mu\leq\Lambda_0.
\end{equation}
We  represent solution as $u_\mu=v+w$, where $v$ solves the problem
$$
\Delta v=0\;\;\mbox{in $\Omega$ and}\;\;v=u_\mu\;\;\mbox{on $\Gamma$}.
$$
The function $v$ belongs to $H^1(\Omega)$ and due to (\ref{Jul11b}) satisfies the estimate
\begin{equation}\label{Jul11c}
||v||_{L_2(\Omega)}\leq C||u_\mu||_{L_2(\Gamma)}\leq C\mu^{-1/2}.
\end{equation}
Therefore $w\in H^1(\Omega,\Gamma)$ and satisfies $-\Delta w=\Lambda(\mu)u_\mu$ in $\Omega$. Multiplying this equation by $w$ and integrating over $\Omega$, we get
$$
\int_\Omega |\nabla w|^2dx=\Lambda(\mu)\int_\Omega (w+v)wdx.
$$
Using the definition of the constant $\Lambda_0$ we derive from the last identity the following estimate
$$
(\Lambda_0-\Lambda(\mu))\int_\Omega |w|^2dx\leq \Lambda(\mu)||w||_{L_2(\Omega)}||v||_{L_2(\Omega)}.
$$
Therefore
$$
||w||_{L_2(\Omega)}\leq \frac{\Lambda(\mu)}{\Lambda_0-\Lambda(\mu)}||v||_{L_2(\Omega)}
$$
or
$$
||v||_{L_2(\Omega)}+||w||_{L_2(\Omega)}\leq \frac{\Lambda_0}{\Lambda_0-\Lambda(\mu)}||v||_{L_2(\Omega)}.
$$
Since the left-hand side $\geq 1$ by using (\ref{Jul11c}), we obtain
$$
\Lambda_0-\Lambda(\mu)\leq C\mu^{-1/2},
$$
which implies (\ref{Jul11a}).

\begin{example}\label{ex1}
Let $\Omega$ be a strip in ${\mathbb R}^2$ i.e $\Omega=\{\ (x,y): x\in {\mathbb R}, ~~  0<y<L\}\ $. 

Consider the following spectral boundary value problem in the cross-section  $[0,L]$. Find Y such that
\begin{equation}\label{eq:helmholtz:seperation:2}
-Y^{\prime\prime}=\lambda Y
\end{equation}
and
 \begin{equation}\label{eq:helmholtz:boundary:1}
 Y^{\prime}(0)-\mu_0 Y(0)= Y^{\prime}(L)+\mu_1 Y(L)=0 ,
\end{equation}
 where   $ \mu_0,\mu_1 $  are non-negative and $ \mu_0+\mu_1 >0$.

The  first eigenvalue of (\ref{eq:helmholtz:seperation:2}) with homogenous Dirichlet boundary conditions in the cross-section  $[0,L]$ is $ \frac{\pi^2}{L^2} $  and we denote it by $ \lambda_0$. Our aim  is to evaluate the first eigenvalue, $\lambda(\mu)$ of problem (\ref{eq:helmholtz:seperation:2}), (\ref{eq:helmholtz:boundary:1}) and to demonstrate that this eigenvalue is close to $ \lambda_0$.
Eigenvalues of  (\ref{eq:helmholtz:seperation:2}) and  (\ref{eq:helmholtz:boundary:1}) are positive. Multiplying both sides of (\ref{eq:helmholtz:seperation:2}) by  $ Y $, integrating by parts and applying the boundary conditions (\ref{eq:helmholtz:boundary:1}) gives,

\begin{equation}\label{eq:helmholtz:boundary:2}
\begin{split}
\Int_{0}^{L}Y^{\prime}(y)^2+\mu_0 Y^2(0)+\mu_1Y^2(L)= \lambda \Int_{0}^{L}Y^2(y) \,dy
\end{split}
\end{equation}
 Therefore  $\lambda$ must  be positive.

 To evaluate the first eigenvalue $\lambda(\mu)$ of (\ref{eq:helmholtz:seperation:2}) and (\ref{eq:helmholtz:boundary:1}), let $\lambda=\alpha^2$, $\alpha>0$. The general  solution to (\ref{eq:helmholtz:seperation:2}) is

\begin{equation}\label{eq:helmholtz:seperation:13}
Y(y)=A\cos(\alpha y)+B\sin(\alpha y)
\end{equation}

Using (\ref{eq:helmholtz:boundary:1}),  we obtain
\begin{equation}\label{eq:helmholtz:seperation:15}
  \alpha B-\mu_0 A=0
\end{equation}
 and
\begin{equation}\label{eq:helmholtz:seperation:16}
A(\mu_1\cos(\alpha L) -\alpha \sin(\alpha L)  ) +B(\mu_1\sin(\alpha L) +\alpha \cos(\alpha L)  )=0
 \end{equation}

For existence of non-trivial solution and for  $\mu=\mu_0=\mu_1 $\ we have

\begin{equation}\label{eq:helmholtz:seperation:17}
  \cot (\alpha L)= \frac{(\alpha^2-\mu^2)}{2 \alpha\mu}
 \end{equation}

 or

\begin{equation}\label{eq:helmholtz:seperation:18}
  \cot(\beta)= \frac{(\beta^2-\mu^2 L^2)} { 2L\beta \mu}
 \end{equation}
 where $\beta=\alpha L $.
 Let
 \begin{equation}\label{eq:helmholtz:seperation:19}
f(\beta) = \frac{(\beta^2-\mu^2 L^2)} { 2L\beta \mu}
  \end{equation}
  The smallest root of (\ref{eq:helmholtz:seperation:18}) is located in the interval  $(0,\pi) $, see Figure \ref{fig:root1} . We denote this root by $\beta_1 $.
   If $ \mu L <\frac{\pi}{2}$, then  $\beta_1<\frac{\pi}{2}$, see Figure \ref{fig:root1} . If $ \mu L >\frac{\pi}{2}$, then  $\beta_1>\frac{\pi}{2}$, see Figure \ref{fig:root2} and $\beta_1 $ tends to  $\pi$  as $\mu$ tends to infinity.

 We denote the first root of (\ref{eq:helmholtz:seperation:17}) by $\alpha_1=\frac{\beta_1}{L}$. Therefore the first eigenvalue  $\lambda(\mu)< (\frac{\pi}{2L})^2 $ when $\alpha_1 > \mu$,  and $\lambda(\mu)> (\frac{\pi}{2L})^2 $
whenever $\alpha_1 < \mu$.
For large $\mu$, we have,
\begin{equation}\label{eq:helmholtz:seperation:20}
f(\beta,\mu) = \frac{\beta}{2L\mu} - \frac{ \mu L}{2\beta}
  \end{equation}
  and $f(\beta,\mu)$ tends to  $-\infty $ as  $\mu \to \infty $. The leading term of the root of (\ref{eq:helmholtz:seperation:18}) is $\pi$, so $\beta_1= \pi -\delta$ where $\delta$ is a small number for large $\mu$. Substituting  $\beta_1$ into (\ref{eq:helmholtz:seperation:18}) and using  (\ref{eq:helmholtz:seperation:20}) for $\mu \to \infty $, we have

 \begin{equation*}
 \frac{\cos (\pi -\delta)}{\sin (\pi -\delta)}= -\frac{1}{\delta}\left(1+\mathcal{O} (\delta^2) \right)=-\frac{\mu L}{2\beta}
   \end{equation*}
   and therefore
    \begin{equation*}
 \delta = \frac{2\beta }{\mu L} +\mathcal{O}( \frac{1}{\mu^3})
   \end{equation*}
As a result, we have that
  $\lambda(\mu) \to \lambda_0 $ as   $ \mu \to \infty $
  
The first eigenvalue  of problem (\ref{eq:helmholtz:seperation:2})  with boundary conditions (\ref{eq:helmholtz:boundary:1}) is given by the formula $\lambda(\mu)=\frac{\beta^2}{L^2} $ and is evaluated from 
equation (\ref{eq:helmholtz:seperation:18}). In Table \ref{table:lambda} we present some numerically computed values of the first eigenvalue for some increasing values of $\mu$. By Inequality \ref{I3aaa}  and Theorem \ref{eq:convergence:equations :1} we can obtain exponential decay of the solution at infinity and convergence for the iterations presented in Section \ref{sec:Two auxiliary problems and alternating iterative procedure}, respectively if $k^2<\lambda(\mu)$.

 \begin{table}[ht] 
\centering
\begin{tabular}{ |c | c | c | c |c |c |c |c |}
\hline 
$\mu$ & 2 & 4 & 6 & 8 & 10 & 12 & 14 \\ [0.5ex]
\hline 
$\lambda$ & 10.6 & 15.6 & 22.6 & 26.3 & 29.7 &32.2 & 34.8 \\
\hline 
\end{tabular}
\caption{ This table presents the first eigenvalue $\lambda(\mu)$ of problem (\ref{eq:helmholtz:seperation:2}) with homogenous Robin boundary conditions in the cross-section  $[0,0.4]$. } 
\label{table:lambda}
 \end{table}

\begin{figure}[!t]
\begin{center}
\includegraphics[width=8cm]{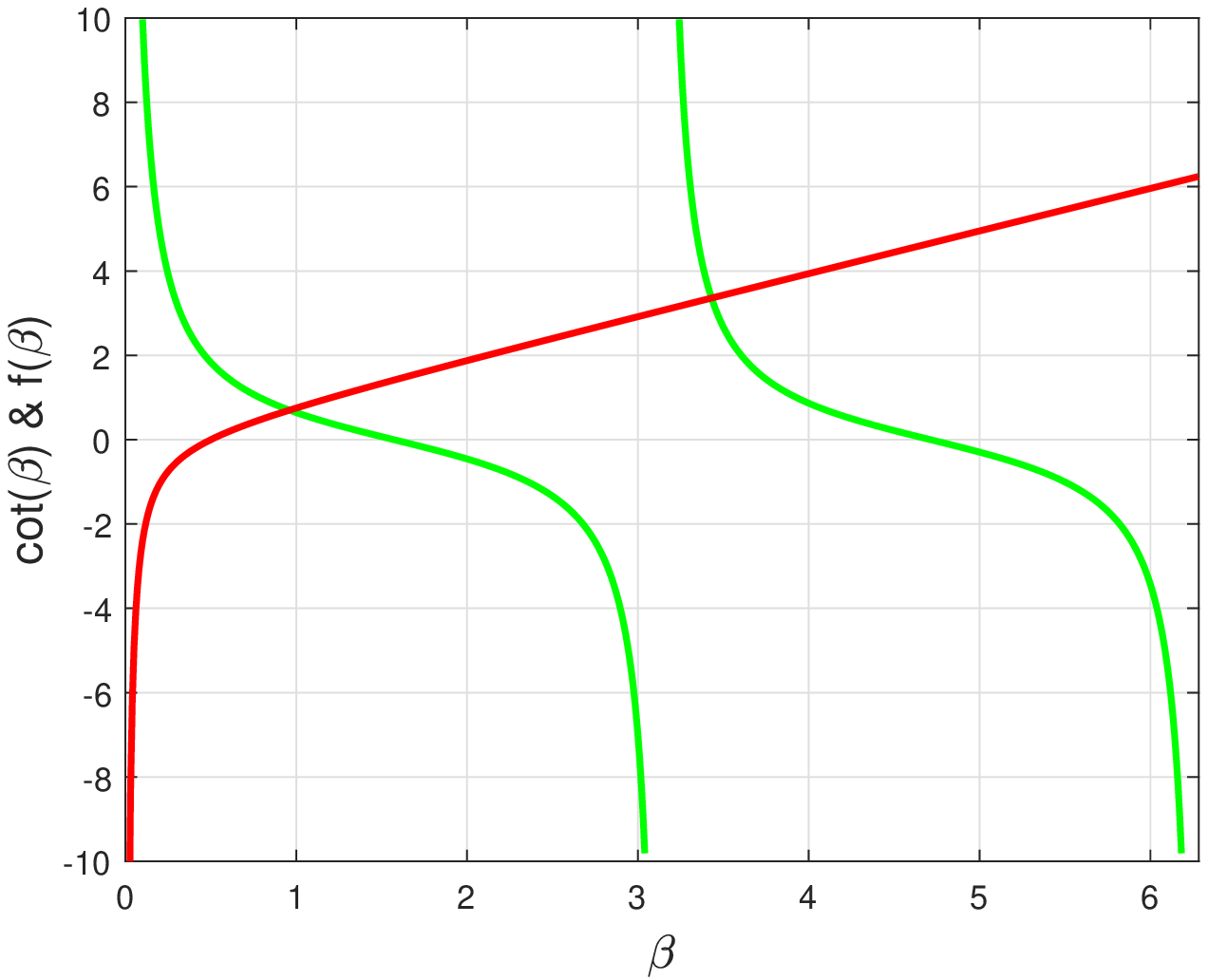}\vspace{10mm}
\end{center}
\caption{\label{fig:root1} Graph of left and right hand side of (\ref{eq:helmholtz:seperation:18}), when
$\mu=3$ and $L=0.5$ which describes the location of the first eigenvalue of problems  (\ref{eq:helmholtz:seperation:2}) and (\ref{eq:helmholtz:boundary:1})  with respect to $\frac{\pi}{2}.$}
\end{figure}
\begin{figure}[!t]
\begin{center}
\includegraphics[width=8cm]{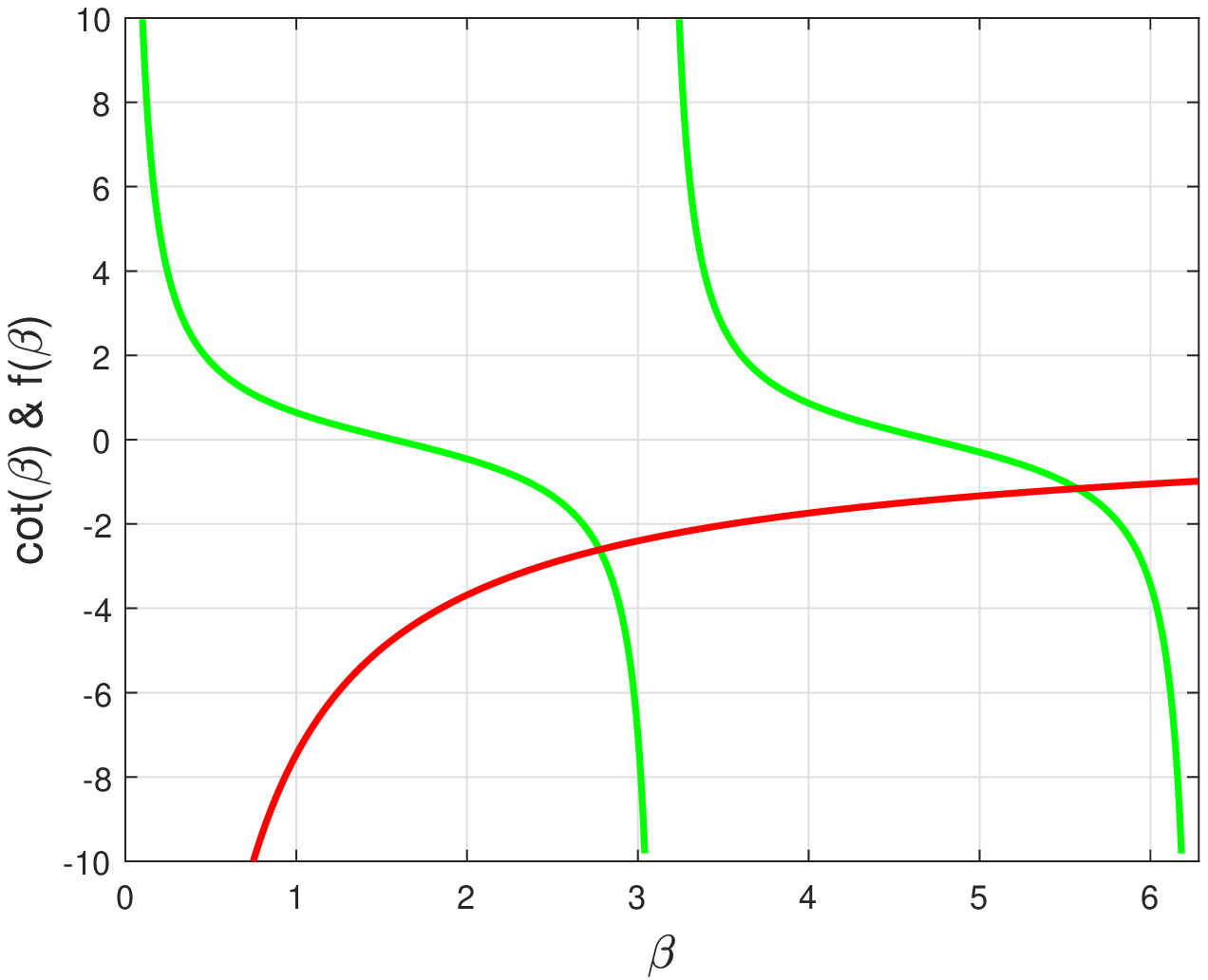}\vspace{10mm}
\end{center}
\caption{\label{fig:root2} Graph of left and right hand side of (\ref{eq:helmholtz:seperation:18}), when
$\mu=30$ and $L=0.5$  which describes the location of the first eigenvalue of problems  (\ref{eq:helmholtz:seperation:2}) and (\ref{eq:helmholtz:boundary:1})  with respect to $\frac{\pi}{2}.$}
\end{figure}

\end {example}
\begin{example}\label{ex2}
Consider the following spectral boundary value problem in  a bounded domain $\omega$ in  $\mathbb R^{d-1}$ with $C^3$ boundary $\partial \omega$.

\begin{equation}\label{eq:spectral:1}
\begin{cases}
-\Delta u_\mu = \lambda_\mu u_\mu & \quad \mbox{in} \quad    \omega \\
 \partial_\nu u_\mu +  \mu u_\mu= 0   & \quad \mbox{on} \quad \partial  \omega  \\
  \end{cases}
 \end{equation}
 where  $\mu\geq 0$,  $\lambda_\mu$ is the least positive eigenvalue and $u_\mu$ is a corresponding eigenfunction.

   \begin{lemma}\label{spectrallemma} The following formula holds
\begin{equation}\label{eq:spectral:2}
 \lambda_\mu=\lambda_0 -\frac{\lambda_1}{\mu}+\mathcal{O}\left( \mu^{-\frac{3}{2}}\right),
 \end{equation}
where   $\lambda_0$  is the least positive eigenvalue of the Dirichlet-Laplacian in $\omega$,
\begin{equation}\label{eq:spectral:6}
 \lambda_1 = \frac{ \int_{\partial \omega}  \mid  \partial_\nu u_0 \mid^2  \,dS}{ \int_{\omega} \mid u_0\mid^2   \,dx},
 \end{equation}
and  $u_0$ is the eigenfunction corresponding to $\lambda_0$.
\end{lemma}
{\bf Proof.}
We normalize the sequence in (\ref{eq:spectral:1}) by $||u_\mu||_{L_2(\omega)}=1$ then from Lemma 3.1 in  \cite{achieng2020}, it follows that
\begin{equation}\label{eq:spectral:3}
\int_\omega |\nabla u_\mu|^2 dx+\mu \int_{\partial \omega}  u_\mu^2 dS \leq C,
\end{equation}
where C does not depend on $\mu$. This implies that $u_\mu $ is weakly convergent  to  $u_D $ in $ H^1(\omega)$, $u_\mu $ is convergent to $u_D $ in $ L^2(\omega)$ and $\lambda_\mu$ converges to $\lambda_D$.

Let us now construct an approximate solution to (\ref{eq:spectral:1}),  $\tilde{u}=u_D+\frac{1}{\mu}u_1$  and  $ \tilde{\lambda}=\lambda_D+\frac{1}{\mu}\lambda_1$, where the function $u_1$ and the number $\lambda_1$ are found from the following  problem
\begin{equation}\label{eq:spectral:4}
\begin{cases}
(\Delta  + \tilde{\lambda}) \tilde{u} =\mu^{-2}\tilde{f} & \quad \mbox{in} \quad    \omega \\
 \partial_\nu \tilde{u}  +  \mu \tilde{u}= \mu^{-1}\tilde{g}  & \quad \mbox{on} \quad  \partial  \omega  \\
  \end{cases}
 \end{equation}
 by equating  coefficients in $\mu^{-1}$ in the equation and in $\mu^0$ in the boundary condition. As the result we get the following equation for  $\lambda_1$ and $u_1$:

 \begin{equation}\label{eq:spectral:5}
\begin{cases}
-\Delta u_1 = \lambda_0 u_1+\lambda_1 u_0 & \quad \mbox{in} \quad \omega \\
u_1+\partial_\nu u_0 =0  & \quad \mbox{on} \quad  \partial  \omega  \\
\end{cases}
\end{equation}
Multiplying the first equation in (\ref{eq:spectral:5}) by $u_0$ and integrating over  $\omega$ we 
obtain the following solvability criterion for the problem (\ref{eq:spectral:5}), which is considered as a problem with respect to $u_1$,
\begin{eqnarray*}
- \int_{\omega} \Delta u_1 u_0 dx &=& \lambda_1  \int_{\omega} u_0^2 dx +\lambda_0 \int_{\omega} u_1u_0 dx\\
&=&\int_{\partial \omega} (-\partial_\nu u_1 u_0 +u_1 \partial_\nu u_0 ) dS -\int_{\omega} u_1 \Delta u_0 dx
\end{eqnarray*}
which implies that
\begin{eqnarray*}
 \lambda_1  \int_{\omega} u_0^2 dx =\int_{\partial \omega} u_1 \partial_\nu  u_0 dS =-\int_{\partial \omega} \mid  \partial_\nu u_0 \mid^2  \,dS
\end{eqnarray*}
or (\ref{eq:spectral:6}).

 The function $u_1$ is determined by solving the problem (\ref{eq:spectral:5}). Clearly, $u_1\in H^1(\omega)$ and
 $$
 \tilde{f}= \lambda_1u_1,\;\;\;\tilde{g}= \partial_\nu u_1.
 $$
Multiplying the first equation in (\ref{eq:spectral:4}) by $u_\mu $, integrating by parts over $\omega $ and applying the boundary conditions in (\ref{eq:spectral:4}) together with (\ref{eq:spectral:1}) we obtain

\begin{eqnarray*}
 \int_{\omega} (\Delta  + \tilde{\lambda}) \tilde{u} u_\mu =\int_{\partial \omega} (\partial_\nu\tilde{u}u_\mu -\tilde{u} \partial_\nu u_\mu ) dS + (\tilde{\lambda}-\lambda_\mu)\int_{\omega} \tilde{u}u_\mu dx &=& \frac{1}{\mu^2} \int_{\omega}\tilde{f} u_\mu dx
\end{eqnarray*}
which implies that
$$(\tilde{\lambda}-\lambda_\mu)\int_{\omega} \tilde{u}u_\mu dx= \frac{1}{\mu^2}\int_{\omega}\tilde{f} u_\mu dx-\frac{1}{\mu}\int_{\partial \omega} \tilde{g} u_\mu dS$$
Since 
$$\int_{\omega} \tilde{u}u_\mu dx=1+\mathcal{O}\left( \frac{1}{\mu} \right)$$
and using (\ref{eq:spectral:3}) we get

$$ \mid \tilde{\lambda}-\lambda_\mu \mid \leq \mathcal{O}( \frac{1}{\mu^2})+C\frac{1}{\sqrt {\mu}}\| \tilde{g}  \|_{L_2(\Omega)}  \leq \frac{C}{\mu^\frac{3}{2}} $$
which proves (\ref{eq:spectral:2}) for $C^3$ boundary  $ \partial  \omega$.
\end {example}

\begin{remark}
We note that inequalities (\ref{I3a}) and (\ref{I3aaa}) implies that solution from $ L^2(\Omega)$ to the Helmholtz equation with Dirichlet  and Robin boundary conditions and with compactly supported right-hand sides exponentially decay at infinity. 
\end{remark}


\section{\label{sec:Solvability of problems (1.5), (1.6)  and convergence of the alternating procedure}Solvability of problems (1.5), (1.6)  and convergence of the alternating iterative procedure}
In this section we describe the function spaces involved in problems (\ref{eq:auxiliary:2}) and (\ref{eq:auxiliary:3}),
define the weak solutions for the problems and state their solvability results.  We also state without proof the theorem on convergence of the alternating iterative procedure described in Section \ref{sec:Two auxiliary problems and alternating iterative procedure}.

\subsection{\label{sec: Function spaces }Function spaces }
The  Sobolev space $H^1(\Omega)$ consists of all functions in $L^2(\Omega)$ whose first order weak derivatives belong to $L^2(\Omega)$.  As an inner product in  $H^1(\Omega)$, we have
\begin{equation}\label{eq:innerproduct:1}
 a_\mu(u,v)= \int_\Omega (\nabla u\cdot \nabla v -k^2 uv)dx +\mu_0 \int_ {\Gamma_0}uv dS +\mu_1 \int_ {\Gamma_1}uv dS,
\end {equation}
The corresponding norm we denote by  $\| u \|_\mu=a_\mu(u,u)^{1/2 }$ and by  Assumption $\ref{I4}$, this norm is equivalent to the standard norm in $H^1(\Omega)$.
We denote by $H^{1/2 }(\Gamma) $, the space of traces of functions  in $H^1(\Omega)$ on  $\Gamma $
. Also,  $H^{1/2 }(\Gamma_0) $, the space of restrictions of functions belonging to $H^{1/2 }(\Gamma) $  to $ \Gamma_0 $  and  $H^{1/2 }_{0}(\Gamma_0) $, the subspace of $H^{1/2 }(\Gamma)$  consisting of functions with supports contained in $ \Gamma_0 $. The dual spaces of $H^{1/2 }_{0}(\Gamma_0) $ is denoted by  $H^{-1/2 }(\Gamma_0)$. Similarly, we can define the spaces  $H^{1/2 }(\Gamma_1) $, $H^{1/2 }_{0 }(\Gamma_1)$ and $H^{-1/2 }(\Gamma_1) $, see  \cite{maxwell2011kozlov, mclean:2000}.

We also define the following subspaces of $H^1(\Omega)$. $ H^{1}(\Omega, \Gamma) $ is the space of functions from $H^1(\Omega)$ vanishing in $\Gamma $. $ H^{1}(\Omega,\Gamma_0 ) $ and $ H^{1}(\Omega,\Gamma_1 ) $ are the spaces of functions from $H^1(\Omega)$ vanishing in $\Gamma _0$ and $\Gamma _1$  respectively.

 \subsection{\label{sec: Weak Solutions} Well-posedness of problems (1.5),(1.6)}
 Similar to \cite{achieng2020}, Section 3.3, we can introduce the weak solution of the Helmholtz equation,
 $ (\Delta+k^2)u=0 $ as a function  $u$ satisfying the following identity
\begin{equation}\label{eq:wellposedness:1}
  \int_\Omega (\nabla u\cdot \nabla v -k^2 uv) dx=0,
\end {equation}
for every function $ v \in  H^{1}(\Omega,\Gamma)$. 

We denote the set of weak solutions to the Helmholtz equation by $\mathbb H$. Clearly $\mathbb H$ is a closed subspace of $H^{1}(\Omega)$. We define the normal derivative  of $u$ as a  function $\partial_\nu u$  in the space $\mathbb H$ satisfying the following inequality,
 
\begin{equation}\label{eq:normal derivative:1}
\| {\partial_\nu u}  \|_{H^{-1/2 }(\Gamma)} \leq C \| u \|_{H^{1}(\Omega)}.
\end {equation}
Consequently, $ \partial_\nu u\Big|_{\Gamma_0}$ and $ \partial_\nu u\Big|_{\Gamma_1}$ are well defined on
$H^{-1/2}(\Gamma_0)$ and $H^{-1/2}(\Gamma_1)$ respectively and satisfy
$$
||\partial_\nu u\Big|_{\Gamma_0}||_{H^{-1/2}(\Gamma_0)}+||\partial_\nu u\Big|_{\Gamma_1}||_{H^{-1/2}(\Gamma_1)}\leq C||\partial_\nu u||_{H^{-1/2}(\Gamma)}.
$$
Again, similar to \cite{achieng2020}, Section 3.4, we define weak solutions for the two boundary value problems (\ref{eq:auxiliary:2}) and (\ref{eq:auxiliary:3}) and show that they are well-posed. We will state the results without proofs since the proofs follow similar arguments as in the case of elliptic equation in bounded domains in  \cite{achieng2020}.
\begin{definition}
\rm The weak solutions for the two boundary value problems (\rm \ref{eq:auxiliary:2}) and (\rm\ref{eq:auxiliary:3}) are defined as follows.
\begin{itemize}
\item[(1)] 
Let $f \in H^{\frac{1}{2}}(\Gamma_0) $ and $\eta\in H^{-1/2 }(\Gamma_1) $.  A function $ u\in H^1(\Omega) $ is a weak solution to (\ref{eq:auxiliary:2}) if
\begin{equation}\label{eq:weaksolution:1}
a_0(u,v) =\int_{\Gamma_1} \eta vds,
\end{equation}
for every function $v\in H^{1}(\Omega, \Gamma_0)$, $u=f$ on $\Gamma_0$ and   $a_0(u,v) $  denote the restrictions of  $a_\mu(u,v) $ to $ H^{1}(\Omega,\Gamma_0 )$.
\item[(2)] Let $ \phi\in H^{\frac{1}{2}}(\Gamma_1)$  and $g \in H^{-1/2 }(\Gamma_0)$.  A function $ u\in H^1(\Omega) $ is a weak solution to (\ref{eq:auxiliary:3}) if
\begin{equation}\label{eq:weaksolution:2}
a_1(u,v)  =\int_{\Gamma_0} g vds,
\end{equation}
for every function $v\in H^{1}(\Omega, \Gamma_1)$, $u= \phi$ on $\Gamma_1$ and $a_1(u,v) $  denote the restrictions of  $a_\mu(u,v) $ to $ H^{1}(\Omega,\Gamma_1 )$
\end{itemize}
\end{definition}
 The solvability results are presented in the following proposition.

\begin{proposition}
\rm Uniqueness conditions for the solutions to the two boundary value problems (\rm\ref{eq:auxiliary:2}) and (\rm \ref{eq:auxiliary:3}) are defined as follows.

\begin{itemize}
\item[(1)] 
Let $f\in H^{1/2 }(\Gamma_0)$  and $\eta\in H^{-1/2 }(\Gamma_1) $ and assume that condition (\ref{I4}) holds for $u \in H^1(\Omega)$ then there exist  a unique weak solution $u \in H^1(\Omega)$ to  problem (\ref{eq:auxiliary:2} ) such that
\begin{equation}\label{eq:well-posed:1}
\| u \|_{H^{1}(\Omega)}\leq C\left(\|f\|_{H^{1/2 }(\Gamma_0)}+\|\eta \|_{H^{-1/2 }(\Gamma_1)}\right),
\end{equation}
where the constant  $C$  is independent of  $f$ and $\eta$.
\item[(2)] 
Let $g \in H^{-1/2 }(\Gamma_0) $ and $ \phi\in H^{\frac{1}{2}}(\Gamma_1)$ and assume that condition (\ref{I4}) holds for $u \in H^1(\Omega)$ then there exist  a unique weak solution $u \in H^1(\Omega)$ to  problem (\ref{eq:auxiliary:3} ) such that
\begin{equation}\label{eq:well-posed:2}
\| u \|_{H^{1}(\Omega)}\leq C\left(\|g\|_{H^{-1/2 }(\Gamma_0)}+\|\phi \|_{H^{1/2 }(\Gamma_1)}\right),
\end{equation}
where the constant  $C$  is independent of  $g$ and $\phi$.
\end{itemize}
\end{proposition}

\subsection{\label{sec: Convergence of the alternating procedure}Convergence of the alternating iterative procedure}

The alternating  iterative algorithm described in section \ref{sec:Two auxiliary problems and alternating iterative procedure} is  linearly dependent on the functions $f $,  $g $ and $\eta $.  This alternating iterative procedure gives a convergent approximation of $u$ in $H^1(\Omega)$ as stated in the following theorem provided that in unbounded domains, condition (\ref{I4}) holds.
\begin{theorem}\label{eq:convergence:equations :1}
Let $f_0\in H^{1/2}(\Gamma_0)$ and $g_0\in H^{-1/2 }(\Gamma_0)$, and let $u\in H^{1}(\Omega)$ be the
solution to the problems {\rm (\ref{I1})},  {\rm (\ref{I2})}. Then for any $\eta\in H^{-1/2}(\Gamma_1)$, the
sequence $\{u_n\}_{n=0}^\infty$, obtained using the algorithm described
in section  \ref{sec:Two auxiliary problems and alternating iterative procedure}, converges to $u$ in $ H^{1}(\Omega)$.
\end{theorem}
The proof of this theorem follows the same lines as the proof of convergence of the alternating iterative procedure in bounded domains presented in Section 4 in \cite{achieng2020} since solutions defined in unbounded domains converges to solutions defined in a bounded cross-section of the unbounded domain.

\section{\label{sec:Numerical Experiments}Numerical Experiments}
In this section, we present numerical experiments and results that illustrate the convergence of the alternating iterative procedure presented in Section \ref{sec:Two auxiliary problems and alternating iterative procedure}.
We specify the geometry and implement a finite difference method to solve the two well-posed boundary problems presented in Section 1.2.

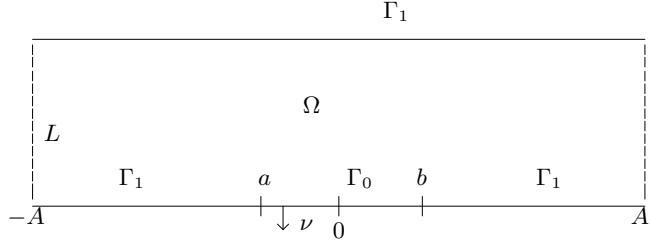
\begin{figure}[!t]
\begin{center}
\setlength{\unitlength}{0.037mm}
\begin{picture}(7253,1605)(0,-10)
\path(390,765)(2590,765)
\path(390,190)(390,245)
\path(390,255)(390,290)
\path(390,300)(390,335)
\path(390,345)(390,380)
\path(390,390)(390,425)
\path(390,435)(390,470)
\path(390,480)(390,515)
\path(390,525)(390,560)
\path(390,570)(390,605)
\path(390,615)(390,650)
\path(390,665)(390,700)
\path(390,710)(390,745)

\path(2590,165)(390,165)
\path(2590,190)(2590,245)
\path(2590,255)(2590,290)
\path(2590,300)(2590,335)
\path(2590,345)(2590,380)
\path(2590,390)(2590,425)
\path(2590,435)(2590,470)
\path(2590,480)(2590,515)
\path(2590,525)(2590,560)
\path(2590,570)(2590,605)
\path(2590,615)(2590,650)
\path(2590,665)(2590,700)
\path(2590,710)(2590,745)
\path(1490,200)(1490,130)
\path(1210,200)(1210,130)
\path(1790,200)(1790,130)
\path(1290,165)(1290,80)
\path(1270,100)(1290,80)
\path(1310,100)(1290,80)
{\small
\put(300,100){\makebox(0,0)[lb]{\smash{$-A$}}}
\put(2540,100){\makebox(0,0)[lb]{\smash{$A$}}}
\put(700,240){\makebox(0,0)[lb]{\smash{$\Gamma_1$}}}
\put(430,400){\makebox(0,0)[lb]{\smash{$L$}}}
\put(1200,240){\makebox(0,0)[lb]{\smash{$a$}}}
\put(1470,50){\makebox(0,0)[lb]{\smash{$0$}}}
\put(2200,240){\makebox(0,0)[lb]{\smash{$\Gamma_1$}}}
\put(1770,240){\makebox(0,0)[lb]{\smash{$b$}}}
\put(1520,240){\makebox(0,0)[lb]{\smash{$\Gamma_0$}}}
\put(1650,840){\makebox(0,0)[lb]{\smash{$\Gamma_1$}}}
\put(1362,500){\makebox(0,0)[lb]{\smash{$\Omega$}}}
\put(1350,80){\makebox(0,0)[lb]{\smash{$\nu$}}}
}
\end{picture}
\caption{\label{fig:Testdomain} Description of the domain considered in the test problem}.
\end{center}
\end{figure}

 For the test we choose a simple rectangular domain that represents an infinite strip truncated at the point where the solution has its support. Let $a, b, A$ and $L$ be real scalars and consider the domain, 

\begin{equation*}
\Omega=(-A,A) \times (0,L) 
\end{equation*}
with 
\begin{equation*}
\Gamma_0=(a,b) \! \times \! \{0\}  \! \quad \mbox{and} \quad \Gamma_1=\left ( (-A,a)\times\{0\} \right)\cup  \! \left ( (b,A) \times \{0\} \right) \! \cup \! \left ( (-A,A) \times \{L\}\right)  \!,
\end{equation*}
see Figure \ref{fig:Testdomain}.

For the test we consider the Cauchy problem for the Helmholtz equation in $\Omega$ as,

\begin{equation}\label{Eq:TestProblem:0}
\begin{cases}
\Delta u(x,y) + k^2 u(x,y) = 0,&\quad \quad -A < x < A, 0< y <L,\\
u(x,0) = f_0(x),& \quad  \quad a\leq x \leq b,\\
u_{y}(x,0) = g_0(x),& \quad \quad a\leq x \leq b,\\
u(-A,y) = u(A,y)=0& \quad  \quad 0\leq y \leq L.
\end{cases}
\end{equation}
 We choose $A$ in such a way that the solution is supported in $(-A, A )\times (0,L) $ and decay exponentially. 
 In the finite difference implementation we introduce a uniform grid on
the domain $\Omega$ of size $N\times M$, such that the step size
is $h=2AN^{-1}$, and thus $M=\textrm{round}(Lh^{-1})$, and use a standard
$\mathcal{O}(h^2)$ accurate finite difference scheme. Our finite difference code solves the Helmholtz equation in the domain $\Omega'=(0,1) \times (0,L')$ and thus we apply the change of variable $x=2A(x'-\frac{1}{2})$  and  $y= 2Ay'$ before computing the numerical solution. Note that the change of variable alters the frequency $k^2$ in the Helmholtz equation and also the Robin boundary conditions since the robin parameters $\mu_0$ and $\mu_1$  are also altered. After having solved the problem on the domain $\Omega'$, we undo the change of variables and display the results in the original domain $\Omega$.

\begin{figure}[!t]
\begin{center}
\includegraphics[width=6.0cm]{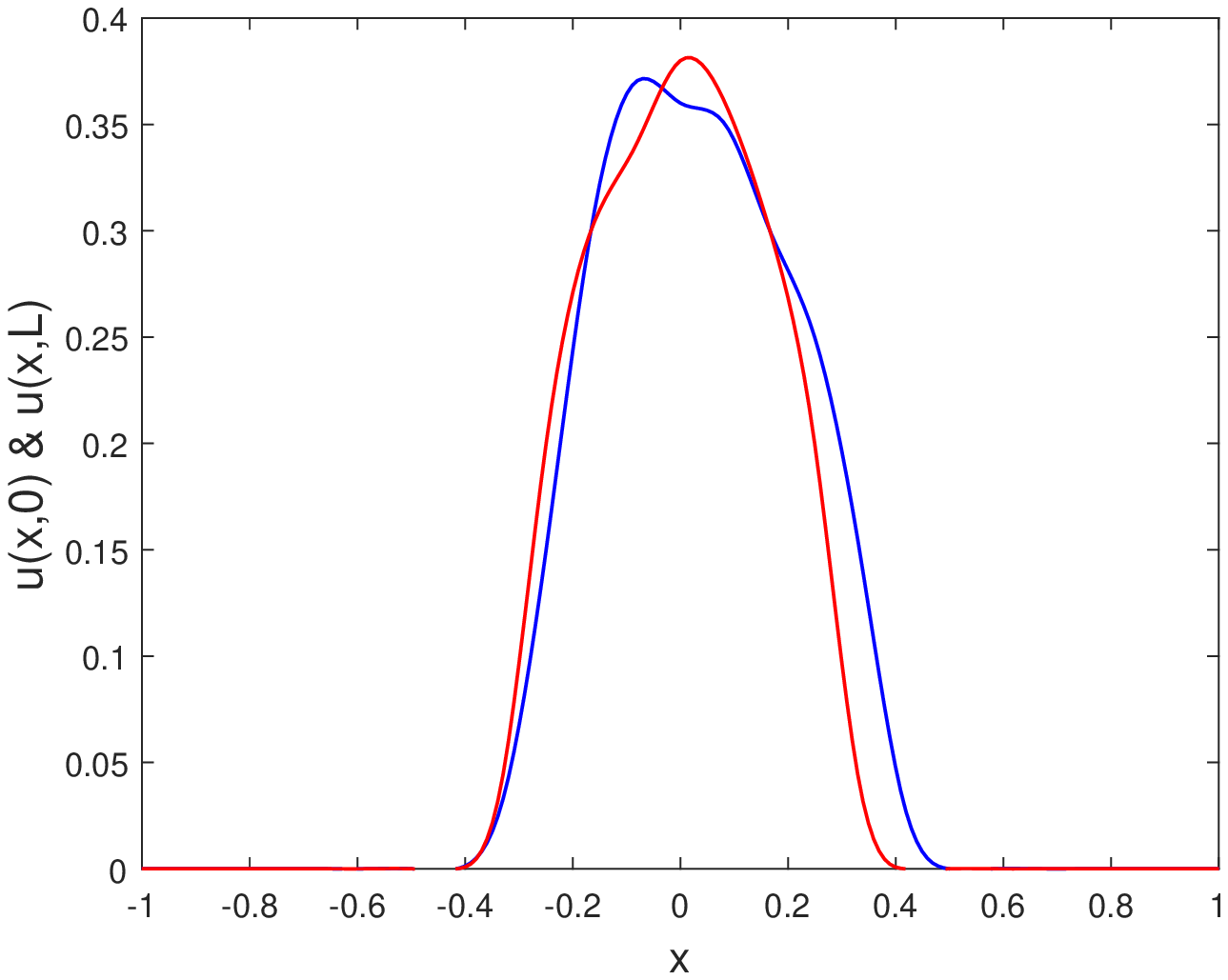}\vspace{3mm}
\includegraphics[width=6.0cm]{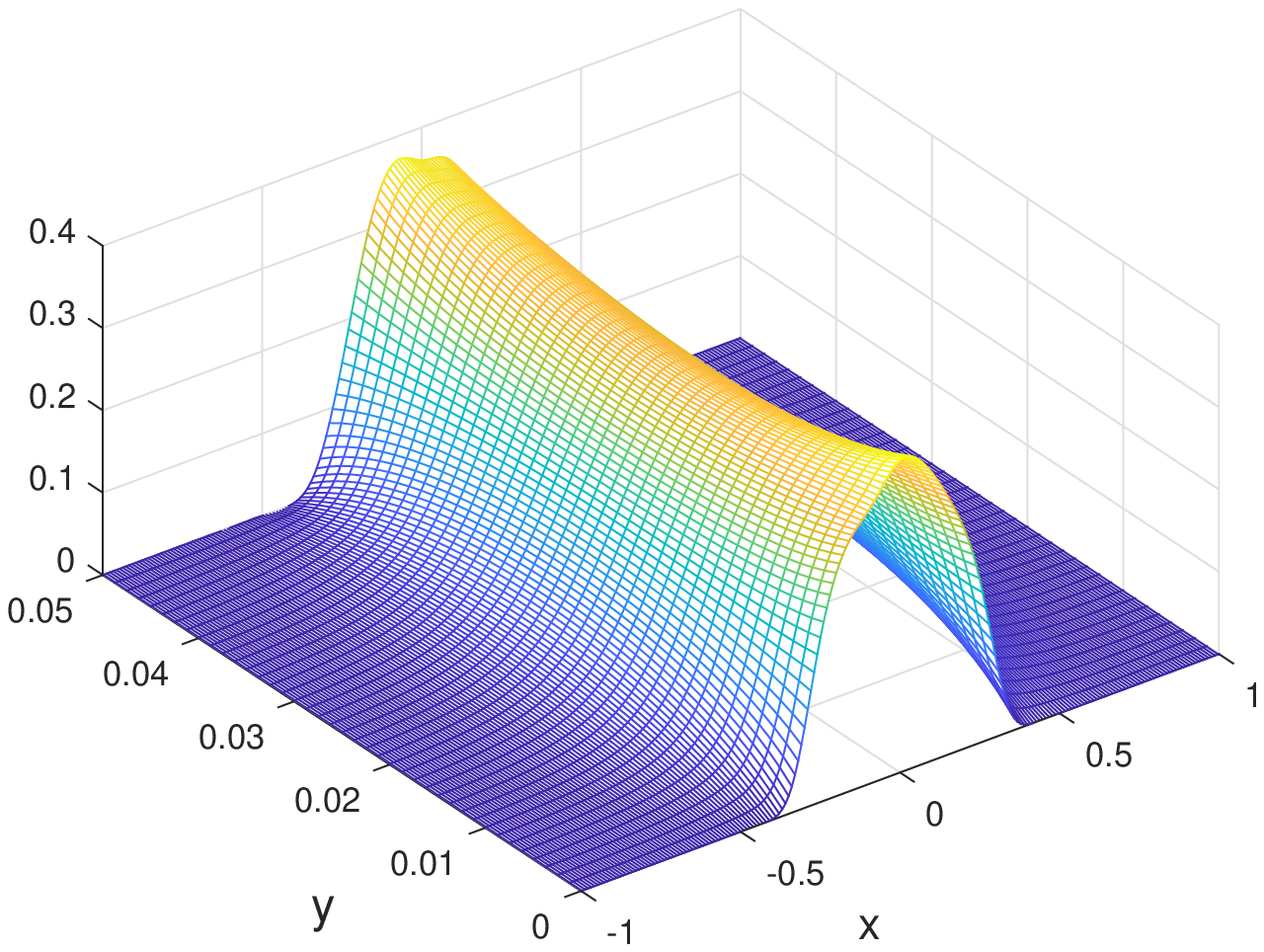}\vspace{3mm}
\end{center}
\caption{\label{fig: testdata} Data $u(x,0)$ (left, red curve),  $u(x,L)$ (left, blue curve) and the numerical solution $u(x,y)$ (right) used to set up the problem.  The graphs are plotted for $-1 \leq x  \leq 1$.}
\end{figure}

In order to obtain the test problem used in our examples, we first set $N=1601$ and $A=4$. We then pick two functions $u(x,0)$  and  $u(x,L)$  which have support in the interval  $[-1,1] $. We use these functions, $A$ and $N$ for all the tests conducted. By solving the Dirichlet problem for the Helmholtz equation in $\Omega$ we obtain a function $u(x,y)$ for $(x,y) \in \Omega$.  In Figure \ref{fig: testdata} we show both the Dirichlet data $u(x,0)$ on $\Gamma_0$, $u(x,L)$  on $\Gamma_1$ and the computed solution for $k^2=5$ and $L =0.4$. For this computation, $M=41$. We display the result in the interval  $-1 \leq x  \leq1$. 

 To illustrate the  Robin-Dirichlet alternating iterative procedure we choose the initial approximations $\eta^{0}(x)=\phi^{0}(x)=0$  and compute a sequence of approximations $\phi^{k}(x)$ for different values of $k^2$, $\mu_0 $, $\mu_1$ and $L $.  We conduct several tests as in the following examples.

\begin{example}\label{test1}
For the first test, we set  $L = 0.4 $,  $\mu_0=\mu_1=2$. We investigate the convergence of the alternating iterative procedure presented in Section \ref{sec:Two auxiliary problems and alternating iterative procedure} with respect to the wavenumber $k^2$. We also investigate how truncation of the domain affect convergence of the procedure.
We observe that for the domain truncated at $A=2$,  the alternating iterative procedure  produces a convergent sequence for $k^2<13.2$ and divergent sequence for  $k^2 >13.2$ while for $A=4$ the procedure produces a convergent sequence for $k^2<12.8$ and divergent sequence for  $k^2 >12.8$. $A=6$ and $A=8$ have no significant differences compared to  $A=4$, see Table \ref{table:truncation}. This motivates our choice for using $A=4$ in all the tests.

We compare the results obtained in test one by the results computed in Table \ref{table:lambda} of Example \ref{ex1}.  This is because according to Theorem \ref{eq:convergence:equations :1} and by Inequality \ref{I3aaa}, convergence of iterations and exponential decay of solution at infinity  is achieved if $k^2<  \lambda( \mu)$, the first eigenvalue of the Robin-Laplacian. Therefore, in unbounded domains, this estimate should determine how the domain is truncated. 
From Table \ref{table:lambda}, for $L = 0.4 $ and  $\mu=\mu_0=\mu_1=2$,  $\lambda(\mu)=10.6$. The variation of this result from the result obtained in test one is possibly due to the choice of A and other errors. 

We also noticed that  the convergence  is quite slow.  See Figure \ref{fig: convergdiverge} and  Figure \ref{fig: testsolution} for illustration of convergence and divergence of the procedure.
 
 \begin{table}[ht] 
\centering
\begin{tabular}{ |c | c | c | c | c |}
\hline 
$A$ & 2 & 4 & 6 & 8  \\ [0.5ex]
\hline 
$k^2 $& 13.2 & 12.8 & 12.7 & 12.7  \\
\hline 
\end{tabular}
\caption{ This table presents the minimum values of  $k^2$ needed for convergence of the Robin-Dirichlet alternating iterative procedure for the domain truncated at different points, A and for  fixed $L = 0.4 $ and $ \mu= \mu_0=\mu_1=2$. } 
\label{table:truncation}
 \end{table} 

 \begin{figure}[!t]
\begin{center}
\includegraphics[width=6.0cm]{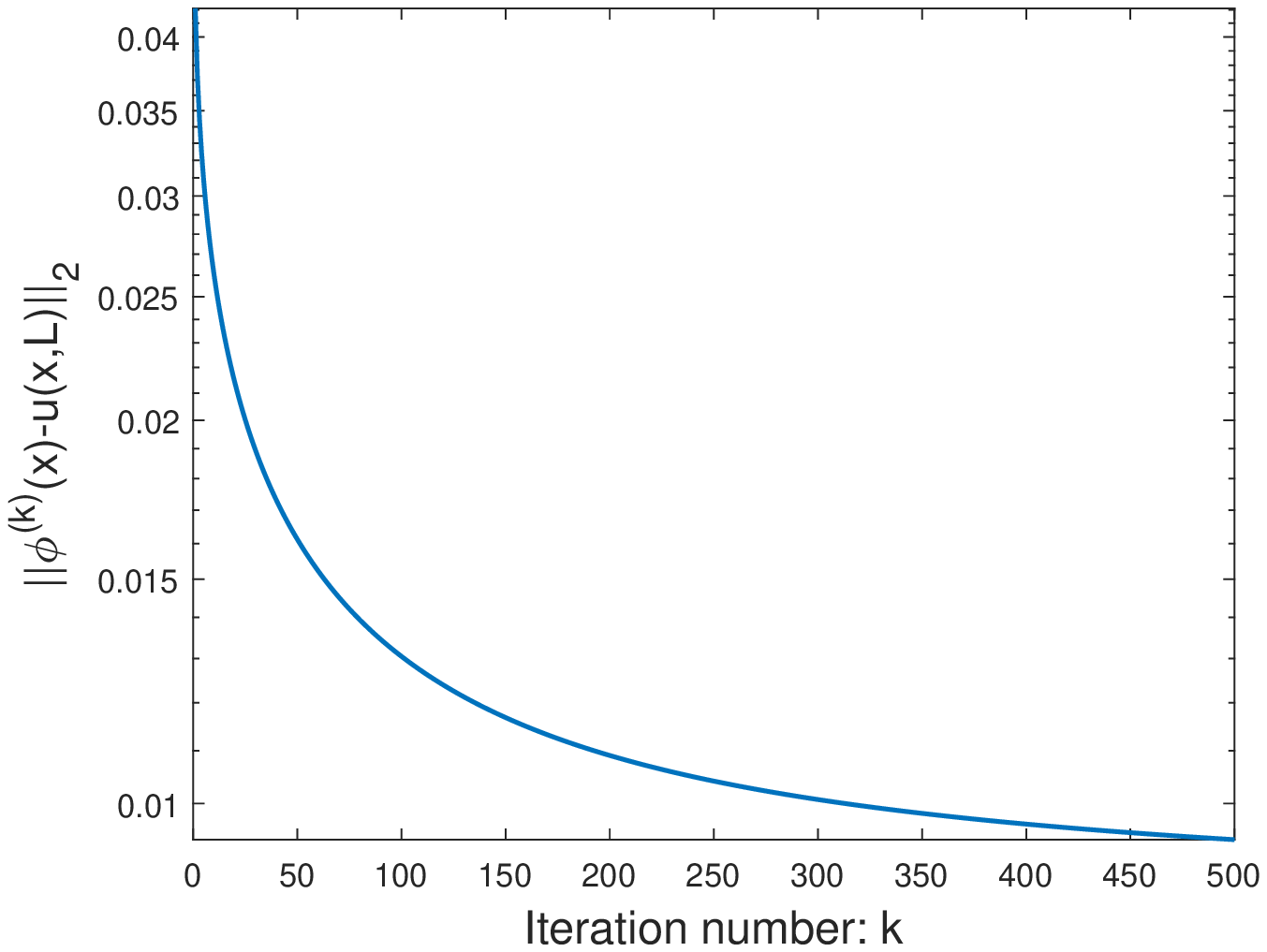}\vspace{3mm}
\includegraphics[width=6.0cm]{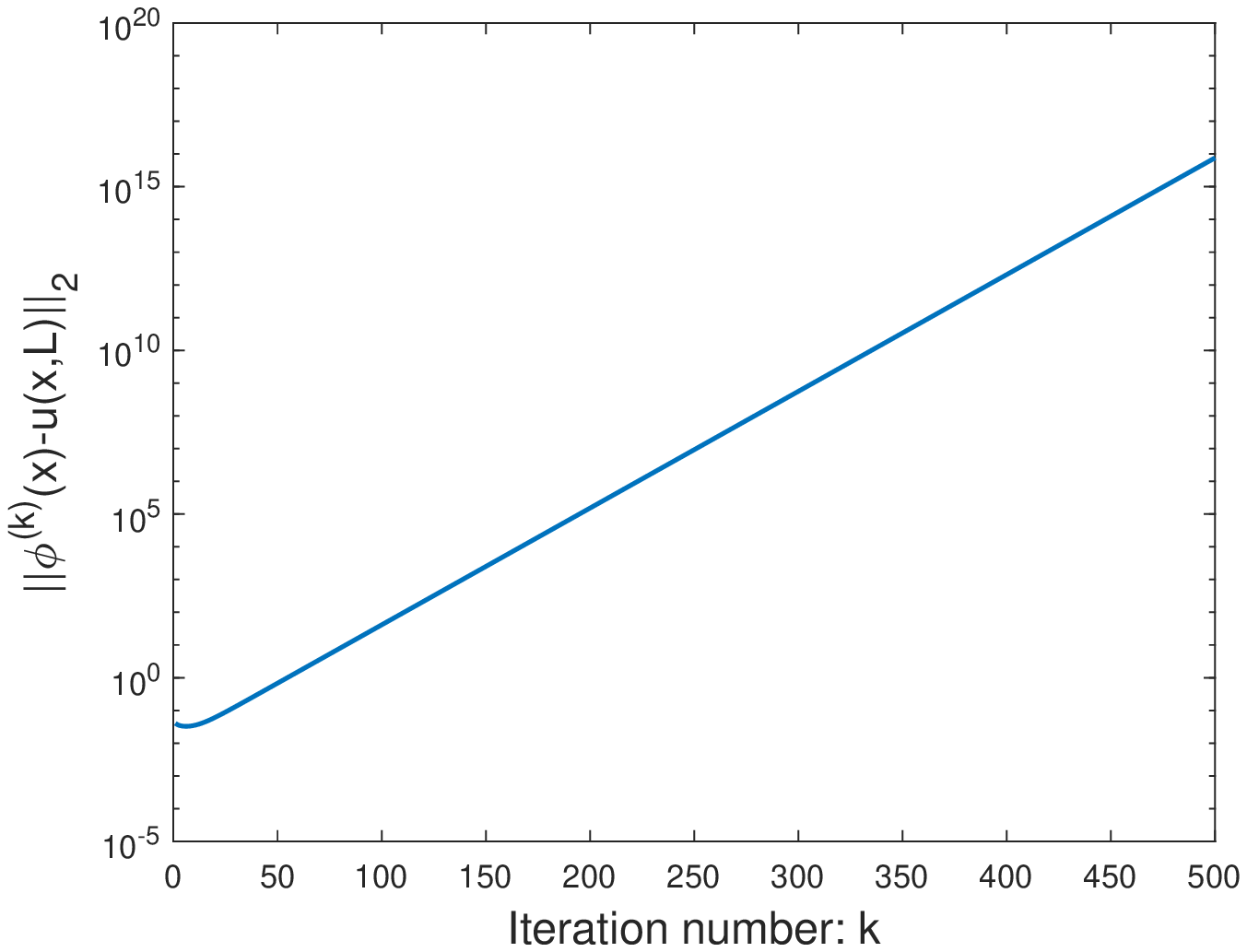}\vspace{3mm}
\end{center}
\caption{\label{fig: convergdiverge}  The error $||\phi^{(k)}-u(x,L)||_2$ during the Robin-Dirichlet iterations for $L =0.4$ and $\mu_0=\mu_1=2.0$. The case $k^2=9.5$ (left) represents a case when the iterations converge and the case $k^2=13.0$ (right) represent when the iterations diverge. The convergent iterations also demonstrate that the convergence is quite slow.}
\end{figure}

 \begin{figure}[!t]
\begin{center}
\includegraphics[width=6.0cm]{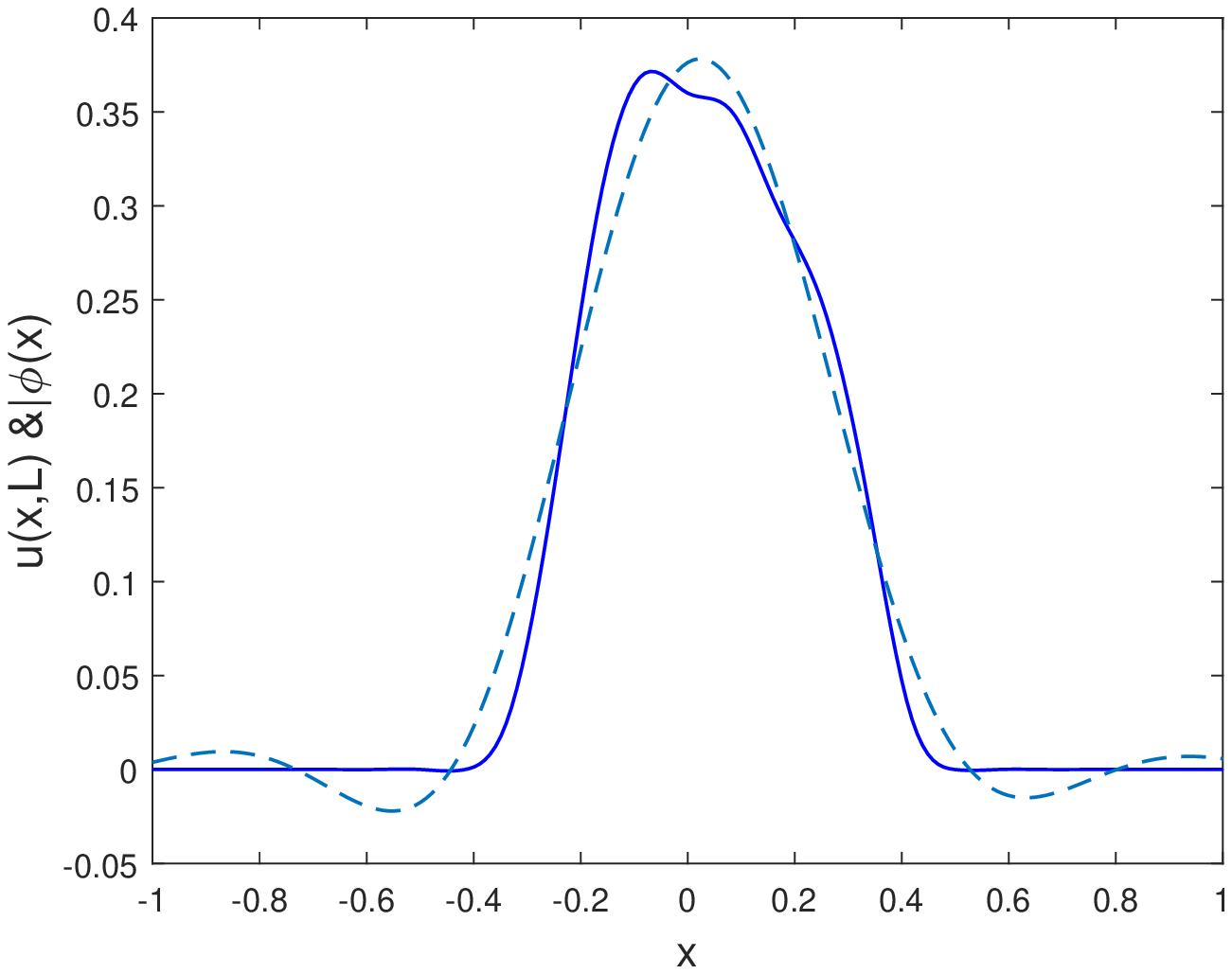}\vspace{3mm}
\includegraphics[width=6.0cm]{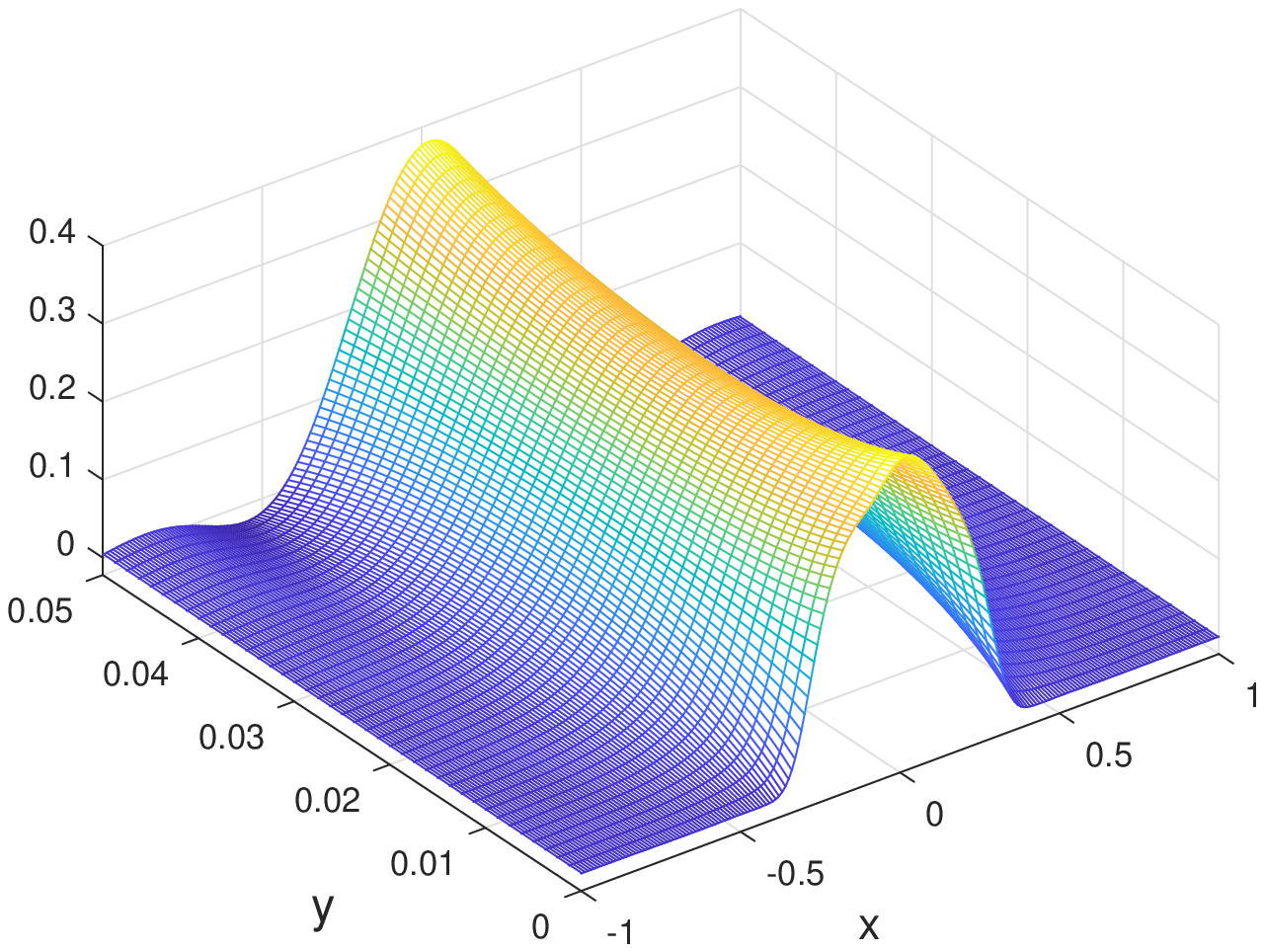}\vspace{3mm}
\end{center}
\caption{\label{fig: testsolution} The exact function $u(x,L)$(left, solid line) and reconstructed function $\phi^j(x)$  for $ j=500 $ (left, dashed line) obtained by the alternating iterative procedure. We also display the approximate numerical solution $u(x,y)$ (right). Here,  $L=0.4 $, $\mu_0=\mu_1=2.0$ and $k^2 =5 $.  We plot the graphs for $-1 \leq x  \leq 1.$} 
\end{figure}
\end {example}
  
\begin{figure}[!t]
\begin{center}
\includegraphics[width=6.0cm]{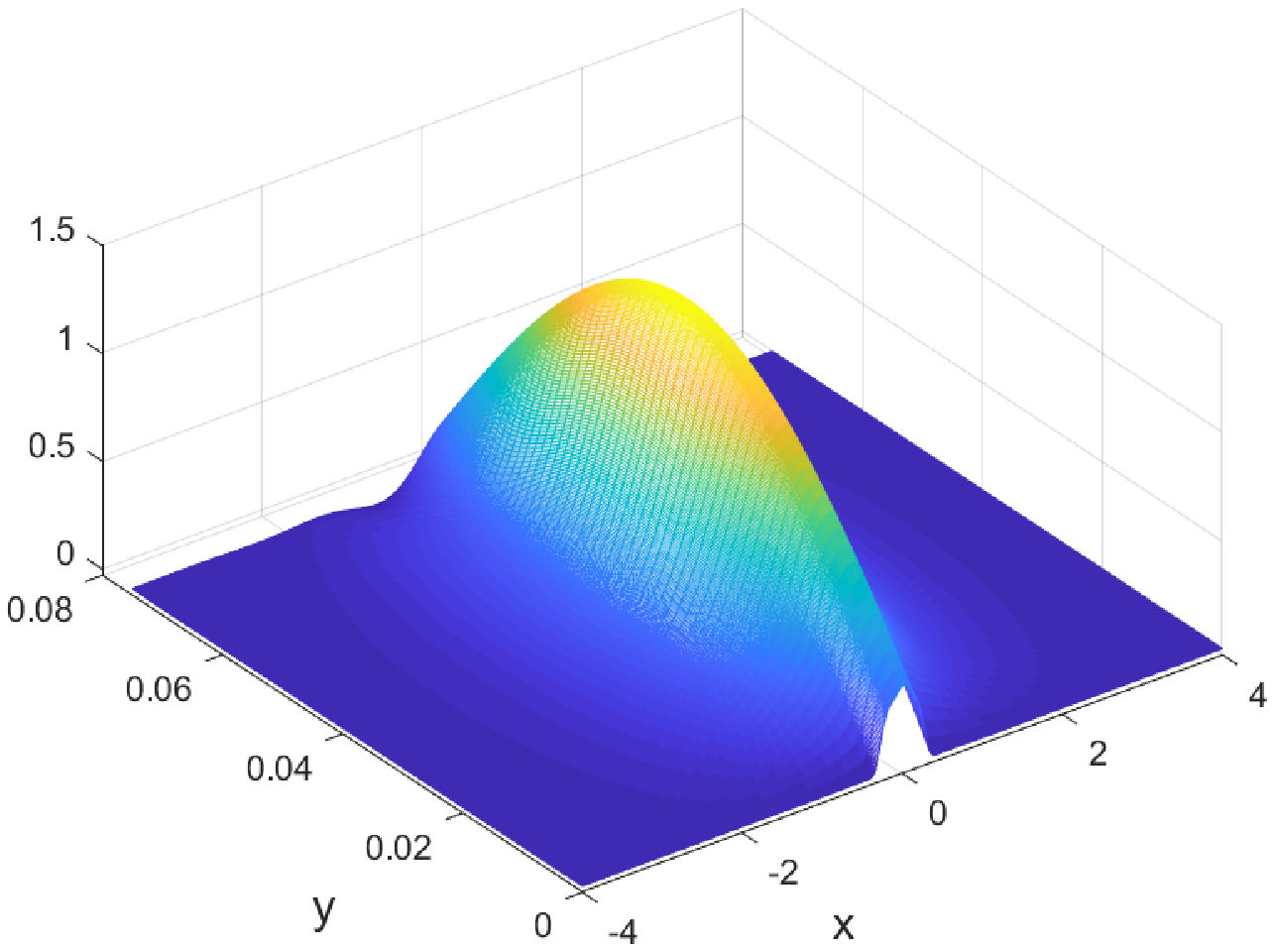}\vspace{3mm}
\includegraphics[width=6.0cm]{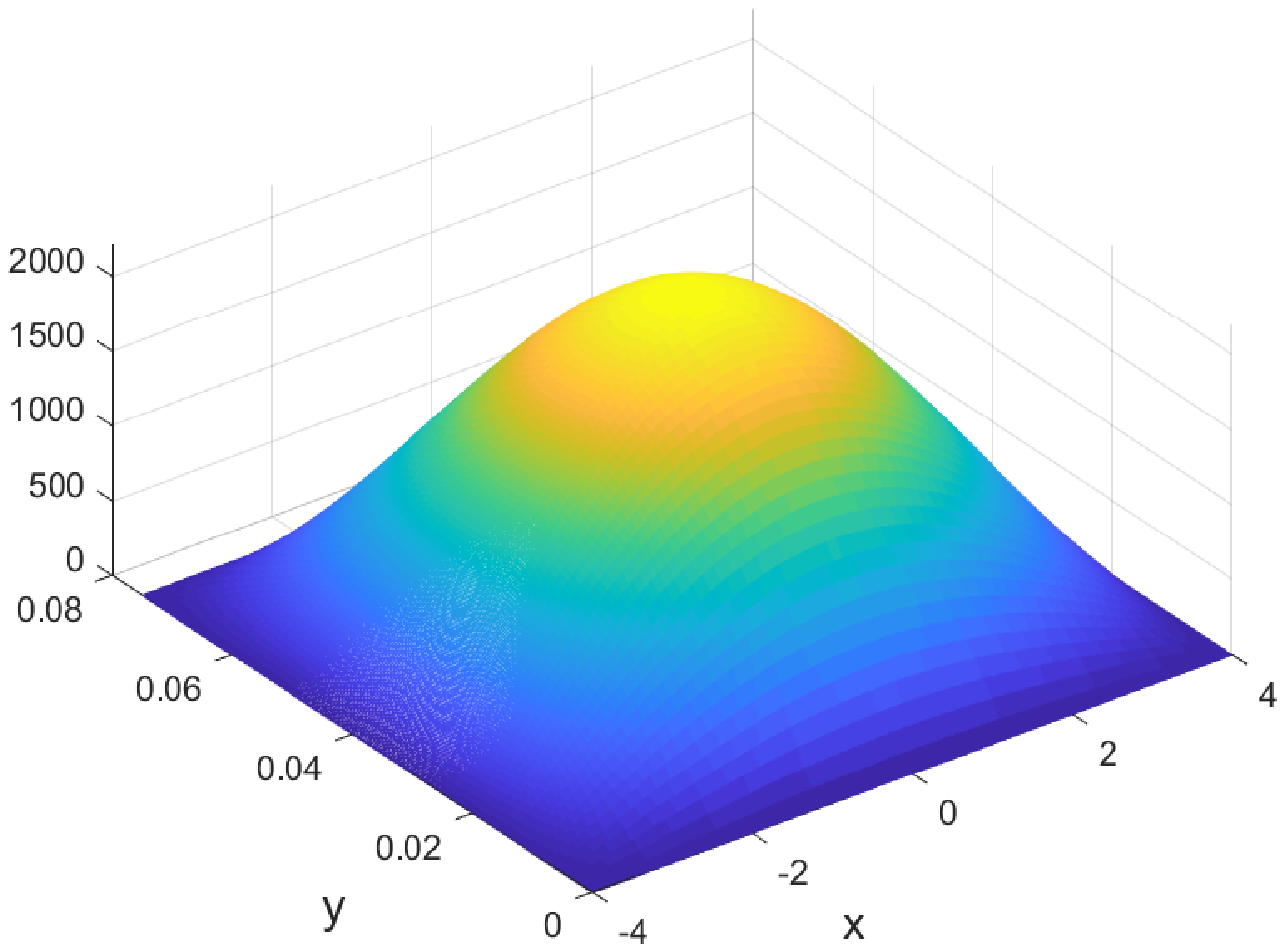}\vspace{3mm}
\includegraphics[width=6.0cm]{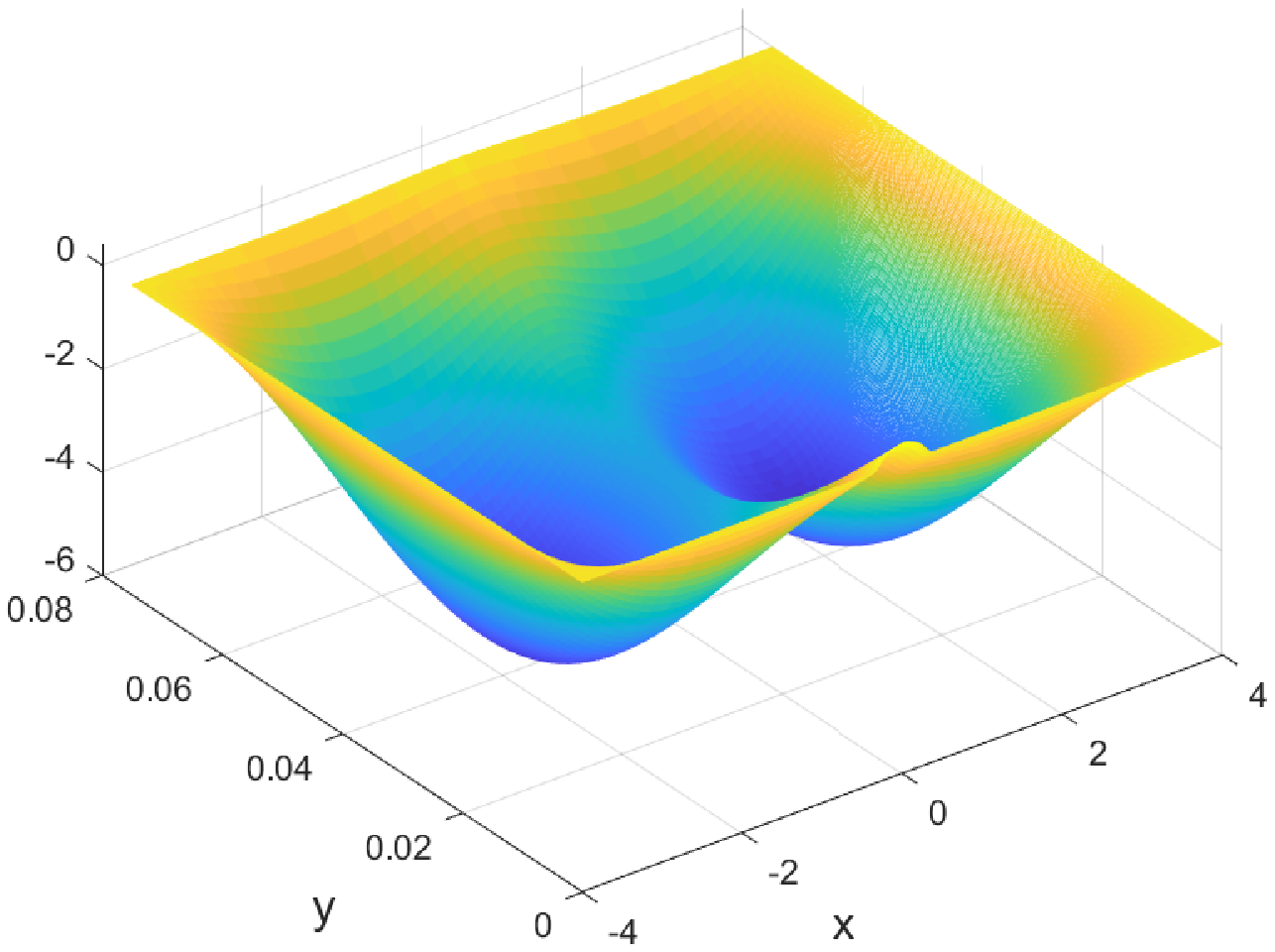}\vspace{3mm}
\includegraphics[width=6.0cm]{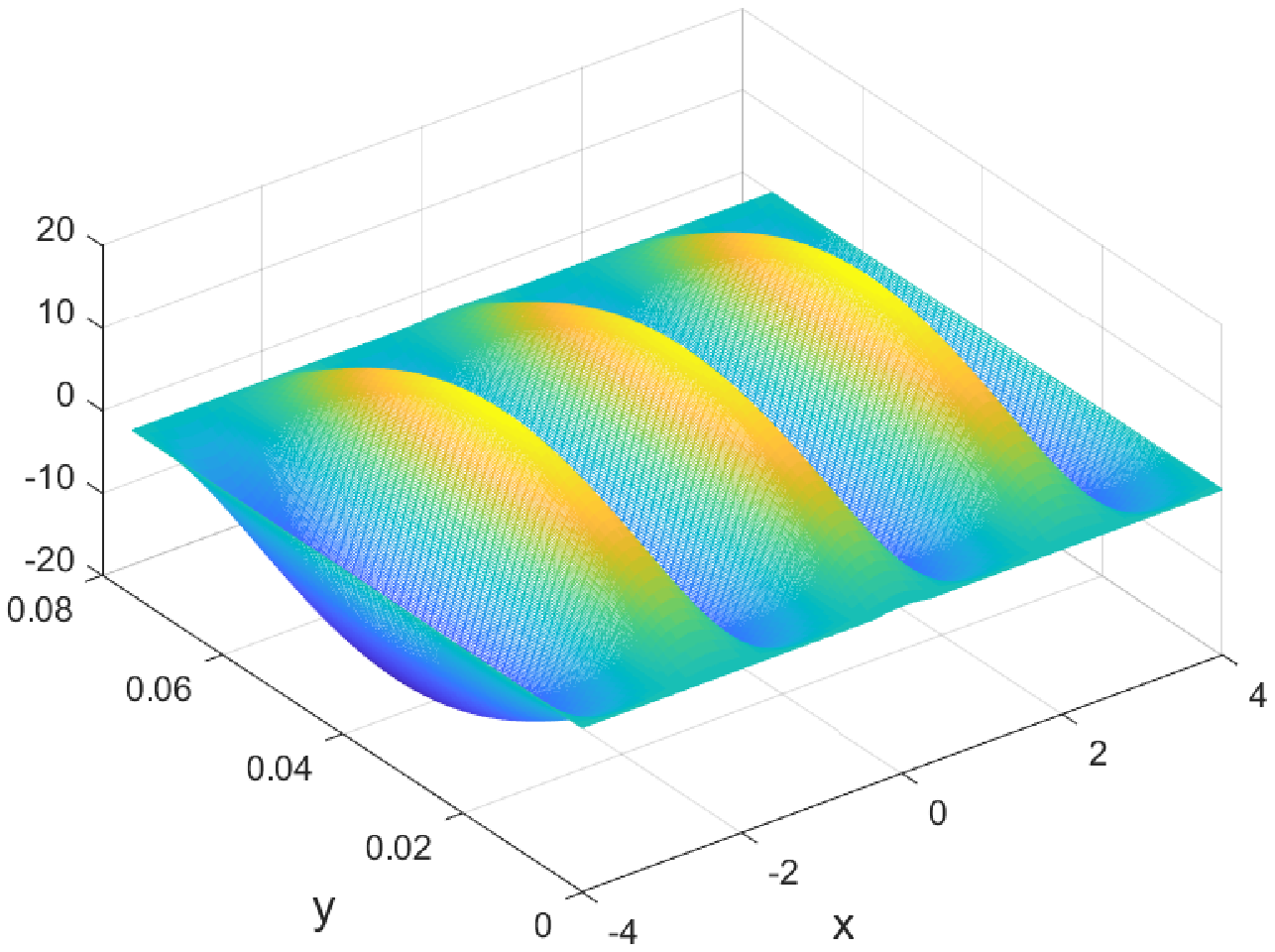}\vspace{3mm}
\end{center}
\caption{\label{fig: badsoln1} Top left is the solution obtained for $k^2=25$ which is supported in $(-1, 1)\times (0,0.6)$ hence solves problem (\ref{Eq:TestProblem:0}). Top right is a sine solution for $k^2=27$ which is just above the limit point where the transition takes place. It is clearly not supported in $(-1, 1)\times (0,0.6)$ hence does not solve problem (\ref{Eq:TestProblem:0}). Bottom left and bottom right are solutions obtained for  $k^2=28$  and  $k^2=35$ respectively which are not supported in $(-1, 1)\times (0,0.6)$ and which illustrate the transition of the solution as $k^2$ increases. All the four solutions are computed for $L=0.6. $ }
\end{figure}

 \begin{example}\label{test2}
 For the second test, we investigate for the minimum $\mu=\mu_0=\mu_1$ needed for convergence of the Robin-Dirichlet alternating iterative procedure for different values of  $L $ and $k^2$.  We set $L=0.2, 0.4$ and $0.6 $ and $k^2$  ranging from $ 5.0$  to $ 50.0$ as shown in Table \ref{table:converg}.
 From the table, we observe that the minimum $\mu$ needed for convergence  increase as $k^2$ increases for all the three different values of $L$. Moreover,
 our problem is ill-posed and the degree of ill-posedness depends on $L $ i.e at  $L =0.2$ we have a better solution compared to  solution obtained at $L =0.6$ which also explains why we need small values of  $\mu$ at $L =0.2$  in order to obtain convergence as compared to the values of $\mu$ at $L =0.6$ needed to obtain convergence.
 
 In Table \ref{table:converg}, $(-)$ means that for $L =0.2$ and  $k^2<15.1$ there is no minimum $\mu$ needed for convergence of the alternating iterative procedure, $\mu=0$ is sufficient. This is as a result of the truncation of the domain so that the Neumann-Dirichlet alternating iterative procedure works for small values of $k^2$ in the Helmholtz equation in a bounded domain. In principle Neumann-Dirichlet alternating iterative procedure should not work in unbounded domains.
 
 Also from Table \ref{table:converg}, $(*)$ represent cases where for large values of $k^2$, we obtain solutions which are not supported in $(-1, 1 )\times (0,L)$ and which do not decay exponentially at infinity hence do not solve the test problem(\ref{Eq:TestProblem:0}).  In Figure \ref{fig: badsoln1} we present four solutions which shows the transition of the solutions as $k^2$ increases.
 
 \begin{table} [ht]
\centering
\begin{tabular}{ |c | c c c |}
\hline
\backslashbox{$ k^2$ }{L}& 0.2 & 0.4 & 0.6 \\ [0.5ex]
\hline 
5.0&-&0.2&1.1 \\ 
10.0&-&1.3 &3.4 \\ 
15.0&-&2.6&7.2 \\ 
20.0&0.5&4.2&15.7 \\
25.0&1.0&6.1 &54.0 \\
30.0& 1.6 & 8.4 & * \\
35.0 &2.1  & 11.9& * \\
40.0 & 2.7 & 16.5 & * \\ 
50.0 & 3.9 & 36.8 & * \\ [1ex]
\hline
\end{tabular}
\caption{ This table presents the minimum $\mu$ needed for convergence of the Robin-Dirichlet alternating iterative procedure for different  values of $L $ and  $k^2.$} 
\label{table:converg}
 \end{table} 
 \end {example}

 \section{\label{sec: Conclusion}Conclusion}
 It was proved in \cite{Lydie:2018} that the Robin-Dirichlet alternating iterative procedure converges even for large values of $k^2$ in the Helmholtz equation in bounded domains if the Robin parameter $\mu$ is appropriately chosen. In this paper, we derive the necessary conditions  for the convergence of the Robin-Dirichlet alternating iterative procedure in unbounded domains. 
 
 Estimate for values of $k^2$ in the Helmholtz equation in terms of positivity of a certain quadratic form which guarantees convergence of the  Robin-Dirichlet alternating iterative procedure is given. We analyse this condition and present some explicit estimates for $k^2$ in terms of eigenvalues of certain auxiliary problems.
 
 For the numerical experiments, we choose a rectangular domain in $\mathbb R^{2}$ that represents a truncated infinite strip and test for convergence of the procedure. 
We know that for the Cauchy problem for the Helmholtz equation in a finite domain, the Neumann-Dirichlet alternating procedure converges  for small values of $k^2$ in the Helmholtz equation while in an infinite domain the Neumann-Dirichlet procedure does not converge at all.
By appropriate truncation of the infinite domain  and with the introduction of  the Robin parameters $\mu_0$ and $\mu_1$, we achieve convergence of the Robin-Dirichlet alternating iterative procedure and the solution decay exponentially at infinity.
However we noticed that the convergence is generally slow. 
 We further investigated dependence of the procedure on the parameter $\mu=\mu_0=\mu_1$ for different values of $k^2$ and $L$ (the distance between the boundaries).
  For $\mu=\mu_0= \mu_1$, the Robin-Dirichlet alternating iterative procedure is divergent for small values of $\mu$ and convergence for large values of $\mu$ and grows as  $k^2$ and $L$  increases. However for large values of $k^2$, we obtain solutions that do not solve the Helmholtz equation in an infinite domain i.e solutions that do not decay exponentially at infinity. 

In our future work, we will investigate the effect of the size of the bounded inclusion on the convergence of the  alternating iterative procedure. We will also investigate the procedure with inexact Cauchy data and add regularization. Finally we will seek to find areas of application.

\bibliographystyle{plain}
\addcontentsline{toc}{anything}{\protect \numberline {References}}
\bibliography{ip,math,helmholtz,ref}
\end{document}